\documentclass{agtart_a}
\pdfoutput=1

\usepackage{amscd,graphicx}

\title[Homology cylinders and the acyclic closure of a free
group]{Homology cylinders and the acyclic closure\\ 
of a free group}

\author{Takuya Sakasai}
\givenname{Takuya}
\surname{Sakasai}
\address{Graduate School of Mathematical Sciences\\
the University of Tokyo\\\newline
3-8-1 Komaba\\Meguro-ku\\Tokyo 153-8914\\Japan}
\email{sakasai@ms.u-tokyo.ac.jp}
\urladdr{}

\volumenumber{6}
\issuenumber{}
\publicationyear{2006}
\papernumber{22}
\startpage{603}
\endpage{631}

\doi{}
\MR{}
\Zbl{}


\received{13 October 2005}
\revised{10 February 2006}
\accepted{23 March 2006}
\published{7 April 2006}
\publishedonline{7 April 2006}
\proposed{}
\seconded{}
\corresponding{}
\editor{}
\version{}

\arxivreference{math.GT/0507260}

\AtBeginDocument{\let\bar\wbar\let\tilde\wtilde\let\hat\what\let\widetilde\wtilde
\let\fulloverline\overline\let\overline\wbar\def\notin{\not\in}}
\newcommand{\we}{\smash{\rlap{\kern 6pt 
\raise 4pt\hbox{\footnotesize $\sim$}}}\longrightarrow}

\makeatletter
\def\cnewtheorem#1[#2]#3{\newtheorem{#1}{#3}[section]
\expandafter\let\csname c@#1\endcsname\c@thm}

\newtheorem{thm}{Theorem}[section]
\cnewtheorem{prop}[thm]{Proposition}
\cnewtheorem{lem}[thm]{Lemma}
\cnewtheorem{cor}[thm]{Corollary}
\newtheorem*{tthm}{\fullref{autsp}}
\theoremstyle{definition}
\cnewtheorem{definition}[thm]{Definition}
\cnewtheorem{example}[thm]{Example}
\cnewtheorem{remark}[thm]{Remark}
\cnewtheorem{question}[thm]{Question}

\makeatother  

\newcommand{\Mg}{\ensuremath{\mathcal{M}_{g,1}}}
\newcommand{\Cg}{\ensuremath{\mathcal{C}_{g,1}}}
\newcommand{\Hg}{\ensuremath{\mathcal{H}_{g,1}}}
\newcommand{\Sg}{\ensuremath{\Sigma_{g,1}}}

\newcommand{\Ker}{\mathop{\mathrm{Ker}}\nolimits}
\newcommand{\Hom}{\mathop{\mathrm{Hom}}\nolimits}
\newcommand{\Image}{\mathop{\mathrm{Image}}\nolimits}
\newcommand{\Id}{\mathop{\mathrm{id}}\nolimits}

\newcommand{\Tmatrix}[1]{\mathop{\left( {#1} \right)}\nolimits}

\newcommand{\Aut}{\mathrm{Aut}}
\newcommand{\End}{\mathrm{End}}
\newcommand{\An}{\mathcal{A}_{n}}
\newcommand{\Bn}{\mathcal{B}_{n}}

\newcommand{\Nil}{F^{\mathrm{nil}}}

\newcommand{\Acy}{F^{\mathrm{acy}}}

\begin{document}

\begin{asciiabstract}
We give a Dehn-Nielsen type theorem for the homology cobordism group
of homology cylinders by considering its action on the acyclic
closure, which was defined by Levine, of a free group. Then we
construct an additive invariant of those homology cylinders which act
on the acyclic closure trivially.  We also describe some tools to
study the automorphism group of the acyclic closure of a free group
generalizing those for the automorphism group of a free group or the
homology cobordism group of homology cylinders.
\end{asciiabstract}

\begin{webabstract}
We give a Dehn--Nielsen type theorem for the homology cobordism group
of homology cylinders by considering its action on the acyclic
closure, which was defined by Levine, of a free group. Then we
construct an additive invariant of those homology cylinders which act
on the acyclic closure trivially.  We also describe some tools to
study the automorphism group of the acyclic closure of a free group
generalizing those for the automorphism group of a free group or the
homology cobordism group of homology cylinders.
\end{webabstract}

\begin{abstract}
We give a Dehn--Nielsen type theorem for the homology cobordism 
group of homology cylinders by considering its action 
on the acyclic closure, 
which was defined by Levine in \cite{le1} and \cite{le2}, 
of a free group. Then we construct an additive invariant of 
those homology cylinders which act on 
the acyclic closure trivially. 
We also describe 
some tools to study the automorphism group of the acyclic 
closure of a free group generalizing 
those for the automorphism group of a free group or the 
homology cobordism group of homology cylinders.
\end{abstract}

\maketitle

\section{Introduction : Dehn--Nielsen's theorem}
\label{ch:intro}
Let $\Sg \, (g \ge 0)$ be a compact connected oriented surface 
of genus $g$ with one boundary component. The fundamental group 
$\pi_1 \Sg$ of $\Sg$ is isomorphic to a free group $F_{2g}$ of 
rank $2g$. We take a word $\zeta \in F_{2g}$ 
which corresponds to the boundary loop of $\Sg$. 

Let $\Mg$ be the mapping class group of $\Sg$ 
relative to the boundary. It is the group of all isotopy classes of 
self-diffeomorphisms of $\Sg$ which fix the boundary pointwise. 
$\Mg$ acts on $\pi_1 \Sg$ naturally, so that we have 
a homomorphism $\sigma\co \Mg \to \Aut F_{2g}$. 
The following theorem due to 
Dehn and Nielsen is well known. 
\begin{thm}[Dehn--Nielsen] The homomorphism $\sigma$ is injective, 
and its image is 
\[\Aut_0 F_{2g} := \left\{ \varphi \in \Aut F_{2g} \mid 
\varphi(\zeta) = \zeta \right\}.\]
\end{thm}
\noindent
From this theorem, we see that an element of $\Mg$ is 
completely characterized by its action on $\pi_1 \Sg$. 
The action of $\Mg$ on $F_{2g}$ induces that on 
its nilpotent quotient 
$N_k:=F_{2g}/(\Gamma^k F_{2g})$ for every $k \ge 2$, 
where $\Gamma^k G$ is the $k^{\rm th}$ term of the lower central series 
of a group $G$ defined by $\Gamma^1 G=G$ and 
$\Gamma^{i} G =[\Gamma^{i-1} G ,G]$ for $i \ge 2$. 
This defines a homomorphism 
$\sigma_k\co \Mg \to \Aut N_k$. 
Note that the restriction of $\sigma_k$ to $\Ker \sigma_{k-1}$, 
whose target is contained in an abelian subgroup of $\Aut N_k$, 
is called the $(k-2)^{\rm nd}$ Johnson homomorphism. 
It has been an important problem 
to determine its image (see \cite{mo}).

Now we consider a generalization of the above argument 
to the homology cobordism group $\Hg$ of 
homology cylinders. The group $\Hg$ has its origin in 
\cite{habi}, \cite{gl}, \cite{le3}, and it 
is regarded as an enlargement of $\Mg$. 
One of the main results of this paper 
is the following Dehn--Nielsen 
type theorem. $\Hg$ has a natural action 
on the group $\Acy_{2g}$ called the acyclic closure of $F_{2g}$, 
which is a completion of $F_{2g}$ in a certain sense defined 
by Levine in \cite{le1} and \cite{le2}, 
and we will determine the image of the representation 
$\sigma^{\mathrm{acy}}\co \Hg \to \Aut \Acy_{2g}$ as follows. 
\begin{tthm}
The image of $\sigma^{\mathrm{acy}}\co \Hg 
\to \Aut \Acy_{2g}$ is 
\[\Aut_0 \Acy_{2g} :=
\{ \varphi \in \Aut \Acy_{2g} \mid
\varphi(\zeta) = \zeta \in \Acy_{2g} \}.\]
\end{tthm}
\noindent
Note that, in this case, 
the homomorphism $\sigma^{\mathrm{acy}}$ 
is not injective. Our next result is the 
construction of a homomorphism 
\[\theta\co \Ker \sigma^{\mathrm{acy}} \longrightarrow 
H_3 (\Acy_{2g})\]
which might be able to detect the elements of 
the kernel of $\sigma^{\mathrm{acy}}$, where the phrase 
``might be'' means that, although we can show that 
this homomorphism is surjective, it is not known, at present, 
whether its target is trivial or not. 
This situation is similar to 
that of some link concordance invariants 
defined by Levine \cite{le1}. 

The group $\Aut \Acy_n$ can be regarded as an enlargement 
of $\Aut F_n$, similar to the situation where 
$\Hg$ enlarges $\Mg$, and it embodies a 
(combinatorial) group-theoretical part of $\Hg$ in the case where $n=2g$. 
From this, we see that $\Aut \Acy_{n}$ itself is an interesting object. 
We will describe some tools to understand this group ---
the Johnson homomorphisms with their refinements, 
and the Magnus representation for $\Aut \Acy_n$. They are 
generalizations of those previously developed 
by Morita \cite{mo} and Kawazumi \cite{ka} for $\Aut F_n$, 
by Habiro \cite{habi}, Garoufalidis--Levine \cite{gl} 
and Levine \cite{le3} for $\Hg$, and by Le Dimet \cite{ld} for 
the Gassner representation of string links. 

The author would like to express his gratitude 
to Professor Shigeyuki Morita 
for his encouragement and helpful suggestions. 
The author also would like to thank the referee for 
useful comments and suggestions. 

This research was partially supported by 
the $21^{\rm st}$-century COE program 
at Graduate School of Mathematical Sciences, 
the University of Tokyo, and by 
JSPS Research Fellowships for Young Scientists.

\section{Definition of homology cylinders}\label{ch:cylinder}
We begin by recalling the definition of 
homology cylinders due to 
Habiro \cite{habi}, Garou\-fal\-idis--Levine \cite{gl} and 
Levine \cite{le3}. 

A {\it homology cylinder} ({\it over} $\Sg$) is 
a compact oriented 3--manifold $M$ equipped 
with two embeddings 
$i_+ , i_-\co  (\Sg,p) \rightarrow (\partial M,p)$ satisfying that 
\begin{enumerate}
\item $i_+$ is 
orientation-preserving and $i_-$ is orientation-reversing, 
\item $\partial M= i_+ (\Sg) \cup i_- (\Sg)$, \ 
$i_+ (\Sg) \cap i_- (\Sg)=i_+ (\partial \Sg)
=i_- (\partial \Sg)$, 
\item $i_+ \bigl|_{\partial \Sg}=i_- \bigl|_{\partial \Sg}$, 
\item $i_+,i_-\co H_{\ast} (\Sg) \rightarrow H_{\ast} (M)$ 
are isomorphisms, 
\end{enumerate}
\noindent
where $p \in \partial \Sg$ is the base point of $\Sg$ and $M$. 
We write a homology cylinder by 
$(M,i_+,i_-)$ or simply by $M$. 
\begin{example}\label{ex:ex1}
$(M,i_+,i_-)=(\Sg \times I, \Id \times 1, \varphi \times 0)$ gives 
a homology cylinder for each $\varphi \in \Mg$, where 
collars of $i_+ (\Sg)$ and $i_- (\Sg)$ are stretched half-way 
along $\partial \Sg \times I$.
\end{example}

Two homology cylinders are said to be {\it isomorphic} if there exists 
an orientation-preserving diffeomorphism between 
the underlying 3--manifolds which is 
compatible with the embeddings of $\Sg$. 
We denote the set of isomorphism classes of homology cylinders 
by $\Cg$. 
Given two homology cylinders $M=(M,i_+,i_-)$ and 
$N=(N,j_+,j_-)$, we can define a new homology cylinder $M \cdot N$ 
by 
\[M \cdot N = (M \cup_{i_- \circ (j_+)^{-1}} N, i_+,j_-).\]
Then $\Cg$ becomes a monoid with the identity element 
\[1_{\Cg}:=(\Sg \times I, \Id \times 1, \Id \times 0).\]
This monoid $\Cg$ is known as an important object to which 
the theory of clasper- or clover-surgeries related to 
finite type invariants of general 3--manifolds is 
applied. 

Instead of the monoid $\Cg$, however, 
we now consider the {\it homology cobordism group $\Hg$ of 
homology cylinders} defined as follows. 
Two homology cylinders $M=(M,i_+,i_-)$ and 
$N=(N,j_+,j_-)$ are {\it homology cobordant} if there exists 
a smooth compact 4--manifold $W$ such that 
\begin{enumerate}
\item $\partial W = M \cup (-N) /(i_+ (x)= j_+(x) , \,
i_- (x)=j_-(x)) \quad x \in \Sg$, 
\item the inclusions $M \hookrightarrow W$, $N \hookrightarrow W$ 
induce isomorphisms on the homology,  
\end{enumerate}
\noindent
where $-N$ is $N$ with opposite orientation. 
Such a manifold $W$ is called a 
{\it homology cobordism} between $M$ and $N$. We denote by $\Hg$ 
the quotient set of $\Cg$ with respect to the equivalence relation 
of homology cobordism. 
The monoid structure of $\Cg$ induces a 
group structure of $\Hg$. In the group $\Hg$, the inverse of 
$(M,i_+,i_-)$ is given by $(-M,i_-,i_+)$. 

The group $\Hg$ has the following remarkable properties. 

First, $\Hg$ contains several groups which relate 
to the theory of low dimensional topology. As we see in 
\fullref{ex:ex1}, we can construct a homology cylinder 
from each element of $\Mg$. This 
correspondence gives an injective 
monoid homomorphism $\Mg \to \Cg$, and moreover, 
the composite of this homomorphism and the natural projection 
$\Cg \to \Hg$ gives an injective group homomorphism. Therefore 
$\Mg$ is contained in $\Hg$. The $g$--component string link 
concordance group $\mathcal{S}_{g}$ is 
also contained in $\Hg$ \cite{le3}. In particular, 
the (smooth) knot concordance group, which coincides with 
$\mathcal{S}_{1}$, 
is contained in $\Hg$. Furthermore, we can inject 
the homology cobordism group $\Theta_{\Z}^3$ of 
homology 3--spheres into 
$\Hg$ by assigning $M \# 1_{\Cg}$ to 
each homology 3--sphere $M$ up to homology cobordism. 

Secondly, we have, as it were, the {\it Milnor--Johnson 
correspondence}, which indicates a similarity between 
the theory of string links and that of homology cylinders. 
Hence we can expect that some methods for studying 
string links (as well as classical knots or links theory in general) 
are applicable to homology cylinders. 

Lastly, we mention about the fundamental group of each homology 
cylinder. For a given homology cylinder $(M,i_+,i_-)$, 
two homomorphisms $i_+,i_-:\pi_1 \Sg \to \pi_1 M$ are not 
generally isomorphisms. 
However, we have the following. 
\begin{thm}[Stallings \cite{st}] \label{thm:st}
Let $A$ and $B$ be groups and $f\co A \rightarrow B$ 
be a $2$--connected homomorphism. Then the induced map 
$f\co A/(\Gamma^k A) \longrightarrow B/(\Gamma^k B)$ 
is an isomorphism for each $k \ge 2$.
\end{thm}
\noindent
Here, a homomorphism $f\co A \to B$ is 
said to be {\it $2$--connected\/} if $f$ induces an isomorphism on 
the first homology, and a surjective homomorphism 
on the second homology. In this paper, the phrase 
``Stallings' theorem'' always means \fullref{thm:st}. 
For each homology cylinder $(M,i_+,i_-)$, 
two homomorphisms 
$i_+,i_-\co F_{2g} \cong \pi_1 \Sg \rightarrow \pi_1 M$ 
are both $2$--connected by definition. 
Hence, they induce isomorphisms on the 
nilpotent quotients of $F_{2g}$ and 
$\pi_1 M$ by Stallings' theorem. 
We write $N_k$ for $F_{2g}/(\Gamma^k F_{2g})$ as before. 
For each $k \ge 2$, 
we define a map $\sigma_k \co  \Cg \rightarrow \Aut N_k$ by 
\[\sigma_k (M,i_+,i_-):= (i_+)^{-1} \circ i_-,\]
which is seen to be a monoid homomorphism. 
It can be also checked 
that $\sigma_k (M,i_+,i_-)$ depends only on 
the homology cobordism class of $(M,i_+,i_-)$, 
so that we have a group homomorphism 
$\sigma_k \co  \Hg \rightarrow \Aut N_k$. 
Note that the restriction of $\sigma_k$ 
to the subgroup $\Mg \subset \Cg$ 
is nothing other than the homomorphism mentioned 
in \fullref{ch:intro}. 
By definition, the image of 
the homomorphism $\sigma_k$ is contained in 
\[\Aut_0 N_k :=
\left\{ \varphi \in \Aut N_k \biggm| 
\begin{array}{c}
\mbox{There exists a lift 
$\widetilde{\varphi} \in \End F_{2g}$ of $\varphi$}\\
\mbox{satisfying $\widetilde{\varphi}(\zeta) \equiv \zeta 
\bmod \Gamma^{k+1}F_{2g}$.} 
\end{array}\right\}. \]
On the other hand, Garoufalidis--Levine 
and Habegger independently showed the following.
\begin{thm}[Garoufalidis--Levine \cite{gl}, Habegger \cite{habe}]
\label{thm:gl} For $k \ge 2$, 
$\Image \sigma_k = \Aut_0 N_k.$ 
\end{thm}
\noindent
Note that the $(k-2)^{\rm nd}$ Johnson homomorphism is obtained by restricting 
$\sigma_k$ to $\Ker \sigma_{k-1}$. From this theorem, 
we can determine its image completely (see \cite{gl}). 

Now we have the following question. Recall that in the case of 
the mapping class group, $\sigma_k\co \Mg \to \Aut N_k$ 
are induced from a single homomorphism 
$\sigma\co \Mg \to \Aut F_{2g}$. 
Then our question is: 
\begin{question}
Does there exist a homomorphism $\Hg \to \Aut \, G$ for 
some group $G$ which induces $\sigma_k\co \Hg \to \Aut N_k$ 
for all $k \ge 2$ ?
\end{question}
\noindent
Some answers to it are given in Remark 2.3 in \cite{gl}. 
Namely, we have a homomorphism 
$\sigma^{\mathrm{nil}}\co \Hg \to \Aut \Nil_{2g}$ 
by combining the homomorphisms $\sigma_k$ for all $k \ge 2$, 
where $\Nil_{2g}:=\displaystyle\varprojlim N_k$ is 
the nilpotent completion of $F_{2g}$. However, $\Nil_{2g}$ 
is too big to treat. Then, the usage of 
the {\it residually nilpotent algebraic closure} of $F_{2g}$, 
which is a countable (as a set) subgroup of $\Nil_{2g}$, 
is also suggested. However, as commented there, we do not know 
whether its second homology is trivial or not. 
The vanishing of it is efficiently used in several situations. In this 
paper, we suggest the usage of 
the {\it acyclic closure $($or HE--closure$)$} 
$\Acy_{2g}$ of $F_{2g}$ to overcome them through 
the argument in subsequent sections.

\section{Construction of an enlargement 
of $\Aut F_n$}\label{ch:enlargement}
Now we fix an integer $n \ge 0$. 
In this section, we define a group, denoted by $\Bn$, 
which can be regarded as an enlargement of $\Aut F_n$. 
The construction of this group 
is analogous to that of the group $\Hg$ of homology cylinders. 
We consider all the arguments in a group level. 
We first construct a monoid $\An$, which enlarges 
the group $\Aut F_n$, and then we obtain the group $\Bn$ by 
taking the quotient of $\An$ by certain equivalence relation. 

\medskip
{\bf Step 1}\qua The construction of the monoid $\An$ 
proceeds as follows. 
Let $\An$ be the set of all equivalence classes of triplets 
$(G,\varphi_+ ,\varphi_- )$ consisting of a finitely presentable 
group $G$ and $2$--connected homomorphisms 
$\varphi_+,\varphi_-\co F_n \rightarrow G$, where 
two triplets $(G,\varphi_+ ,\varphi_- )$, 
$(G',\psi_+ ,\psi_- )$ are said to be 
{\it equivalent} if there exists an isomorphism 
$\rho\co G \stackrel{\cong}{\longrightarrow} G'$ 
which makes the following diagram commutative:
$$\includegraphics{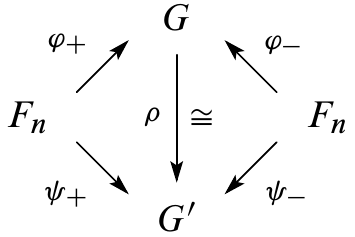}$$
Note that for such a triplet $(G,\varphi_+,\varphi_-)$, 
homomorphisms $\varphi_+$ and $\varphi_-$ are injective, which 
follows from Stallings' theorem and the fact 
that $F_n$ is residually nilpotent. 

We define a multiplication $\cdot$ on $\An$ as follows:
\[\begin{array}{cccc}
 &\An \times \An & \to & \An \\
 & \mbox{\rotatebox[origin=c]{90}{$\in$}} & &
\mbox{\rotatebox[origin=c]{90}{$\in$}}\\
&\left((G, \varphi_+, \varphi_-),(G', \psi_+, \psi_-) \right)
&\mapsto&(G \ast_{F_n} G', \varphi_+, \psi_-)
\end{array}
\]
where $G \ast_{F_n} G'$ is obtained by taking the amalgamated product 
of $G$ and $G'$ with respect to 
$\varphi_-\co F_n \to G$ and $\psi_+\co F_n \to G'$. 
Since $\varphi_-$ and $\psi_+$ are both injective, we can use 
the Mayer--Vietoris exact sequence for the homology of 
amalgamated products (see \cite{br}), so that the above map gives 
a well-defined monoid structure of $\An$ with the identity element 
$(F_n,\Id,\Id)$. 

\begin{example}
$\Aut F_n$ can be seen as a submonoid of $\An$ 
by assigning to each automorphism $\varphi$ of $F_n$ 
an element $(F_n, \Id, \varphi)$ of $\An$. 
This correspondence gives an injective monoid homomorphism 
as shown in \fullref{injective}.
\end{example}

\begin{example}\label{ex}
Consider the monoid of all $2$--connected 
endomorphisms of $F_n$. 
To each $2$--connected endomorphism $\varphi$ of $F_n$, 
we can assign $(F_n, \Id, \varphi) \in \An$. 
This correspondence is also an injective monoid homomorphism. 
For example, consider an endomorphism 
$\psi\co F_2=\langle x_1,x_2 \rangle 
\rightarrow F_2$ defined by 
\[\psi(x_1)=x_1 x_2 x_1 x_2^{-1} x_1^{-1}, \qquad 
\psi(x_2)=x_2 \]
where we take $n=2$. 
As we will see in \fullref{twoconn}, $\psi$ is not 
an automorphism of $F_2$ but a $2$--connected endomorphism. 
Hence $(F_2, \Id, \psi)$ gives an example 
of elements of $\mathcal{A}_2$ 
which are not contained in $\Aut F_2$. 
\end{example}

\begin{example}\label{cgtoa}
For each homology cylinder $(M,i_+,i_-)$, we can obtain 
an element $(\pi_1 M,i_+,i_-)$ of $\mathcal{A}_{2g}$. 
This correspondence gives a monoid homomorphism 
$\Cg \to \mathcal{A}_{2g}$. 
\end{example}

{\bf Step 2}\qua We construct the group $\Bn$ from 
the monoid $\An$ as follows. 
Two elements $(G, \varphi_+, \varphi_-)$, $(G', \psi_+, \psi_-)$ 
of $\An$ are said to be {\it homology cobordant} if 
there exist a finitely presentable group $\widetilde{G}$ and 
$2$--connected homomorphisms
\[\varphi\co G \longrightarrow \widetilde{G}, \qquad 
\psi\co G' \longrightarrow \widetilde{G}\]
which make the following diagram commutative:
$$\includegraphics{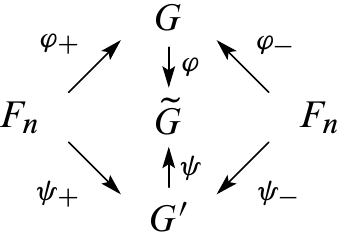}\leqno{\raise 30pt\hbox{($*$)}}$$
We define $\Bn$ to be the quotient set of $\An$ 
with respect to the equivalence relation generated by 
the relation of homology cobordism. 
Then we can endow with a group structure on $\Bn$ 
from the monoid structure of $\An$. 
In $\Bn$, the inverse element of $(G,\varphi_+, \varphi_-)$ is given by 
$(G,\varphi_-, \varphi_+)$. Indeed, 
$G \ast_{F_n} G \xrightarrow{\Id \ast_{F_n} \Id} G 
\xleftarrow{\varphi_+} F_n$ 
gives a homology cobordism between 
$(G \ast_{F_n} G ,\varphi_+, \varphi_+)=(G,\varphi_+, \varphi_-)\cdot
(G,\varphi_-, \varphi_+)$ and $(F_n,\Id,\Id)$. 

\begin{example}\label{hgtob}
The monoid homomorphism $\Cg \to \mathcal{A}_{2g}$ considered 
in \fullref{cgtoa} induces a group homomorphism 
$\Hg \to \mathcal{B}_{2g}$. This homomorphism gives 
an enlargement of the inclusion 
$\sigma\co  \Mg \to \Aut F_{2g}$. 
\end{example}

Fundamental properties of the group $\Bn$ will be 
mentioned in \fullref{ch:properties} after 
seeing a relationship with the acyclic closure of 
a free group.

\section{The acyclic closure of a group}\label{ch:acy}

The concept of the acyclic closure (or HE--closure in \cite{le2}) 
of a group 
was defined as a variation of the algebraic closure of a group 
by Levine in \cite{le1}, \cite{le2}. We briefly summarize the 
definition and fundamental properties. 
We also refer to Hillman's book \cite{hi}. 
The proofs of the propositions in this section 
are almost the same as those for the algebraic closure in \cite{le1} 
(see \fullref{rmk}). 
\begin{definition}
Let $G$ be a group, and let $F_n=\langle x_1, x_2, \ldots, x_n 
\rangle$ be a free group of rank $n$. 

(i)\qua We call each element $w=w(x_1,x_2, \ldots, x_n)$ 
of $G \ast F_n$ a {\it monomial}. A monomial $w$ is 
said to be  {\it acyclic} if 
\[w \in \Ker \left( G \ast F_n \xrightarrow{\mathrm{proj}} F_n 
\longrightarrow H_1 (F_n) \right).\]
(ii)\qua Consider the following ``equation'' with variables 
$x_1,x_2,\ldots,x_n$:
\[
\left\{\begin{array}{ccl}
x_1 & = & w_1(x_1, x_2, \ldots, x_n)\\
x_2 & = & w_2(x_1, x_2, \ldots, x_n)\\
& \vdots & \\
x_n & = & w_n(x_1, x_2, \ldots, x_n)
\end{array}\right. .
\]
When all monomials $w_1,w_2,\ldots,w_n$ are 
acyclic, we call such an equation an {\it acyclic system} 
over $G$.
(iii)\qua A group $G$ is said to be {\it acyclically closed} if 
every acyclic system over $G$ with $n$ variables has 
a unique solution in $G$ for any $n \ge 0$. 
\end{definition}
\noindent
We denote the phrase ``acyclically closed'' by AC, for short. 
\begin{example}
Let $G$ be an abelian group. For $g_1,g_2,g_3 \in G$, 
consider the equation
\[\left\{\begin{array}{l}
x_1=g_1 x_1 g_2 x_2 x_1^{-1} x_2^{-1}\\
x_2=x_1 g_3 x_1^{-1}
\end{array}\right. ,
\]
which is an acyclic system. Then we have a unique solution 
$x_1=g_1 g_2 , \, x_2=g_3$.
\end{example}
\noindent
As we can expect from this example, all abelian groups are AC. 
In fact, all nilpotent groups and the nilpotent completion 
of a group are AC, which are deduced from the following 
fundamental properties of AC--groups and the fact that 
the trivial group is AC. 
\begin{prop}[Proposition 1 in \cite{le1}] 
${\rm (a)}$\qua Let $\{ G_{\alpha} \}$ be a family of 
AC--subgroups of an AC--group $G$. Then $\bigcap_\alpha G_\alpha$ is 
also an AC--subgroup of $G$.

${\rm (b)}$\qua Let $\{ G_{\alpha} \}$ be a family of AC--groups. 
Then $\prod_\alpha G_\alpha$ is also an AC--group.

${\rm (c)}$\qua When $G$ is a central extension of $H$, then 
$G$ is an AC--group if and only if $H$ is an AC--group.

${\rm (d)}$\qua For each direct system $($resp.\ inverse system$)$ 
of AC--groups, the direct limit $($resp.\ inverse limit$)$ is 
also an AC--group.
\end{prop}

Next we define the acyclic closure of a group. 
\begin{prop}[Proposition 3 in \cite{le1}]
For any group $G$, there exists a pair of 
a group $G^{\mathrm{acy}}$ and 
a homomorphism $\iota_G\co G \rightarrow G^{\mathrm{acy}}$ 
satisfying the following properties:
\begin{enumerate}
\item $G^{\mathrm{acy}}$ is an AC--group.
\item Let $f\co G \rightarrow A$ be a homomorphism 
and suppose that $A$ is an AC--group. 
Then there exists a unique homomorphism 
$f^{\mathrm{acy}}\co G^{\mathrm{acy}} \rightarrow A$ 
which satisfies $f^{\mathrm{acy}} \circ \iota_G = f$. 
\end{enumerate}
\noindent
Moreover such a pair is unique up to isomorphisms.
\end{prop}
\begin{definition}
We call $\iota_G$ (or $G^{\mathrm{acy}}$) obtained above 
the {\it acyclic closure} of $G$.
\end{definition}
\noindent
Taking the acyclic closure of a group is functorial, namely, 
for each group homomorphism $f\co G_1 \to G_2$, we obtain 
a homomorphism $f^{\mathrm{acy}}\co G_1^{\mathrm{acy}} \to 
G_2^{\mathrm{acy}}$ by applying the universal property of 
$G_1^{\mathrm{acy}}$ to the homomorphism 
$\iota_{G_2} \circ f$, and the composition of homomorphisms 
induces that of the corresponding homomorphisms 
on acyclic closures. 

The most important properties of the acyclic closure 
are the following. 
\begin{prop}[Proposition 4 in \cite{le1}]\label{2cn}
For every group $G$, the acyclic closure 
$\iota_G\co G \rightarrow G^{\mathrm{acy}}$ is $2$--connected.
\end{prop}
\begin{prop}[Proposition 5 in \cite{le1}]\label{isom}
Let $G$ be a finitely generated group and $H$ be 
a finitely presentable group. For each 
$2$--connected homomorphism 
$f\co G \rightarrow H$, the induced homomorphism 
$f^{\mathrm{acy}}\co G^{\mathrm{acy}} 
\to H^{\mathrm{acy}}$ on 
acyclic closures
is an 
isomorphism. 
\end{prop}
\noindent
From \fullref{2cn} and Stallings' theorem, the 
nilpotent quotients of a group and 
those of its acyclic closure are isomorphic. Note that 
the homomorphism $\iota_G$ is not necessarily injective. 

\begin{prop}[Proposition 6 in \cite{le1}]\label{seq}
For any finitely presentable group $G$, 
there exists a sequence of finitely presentable groups 
and homomorphisms
\[G=P_0 \rightarrow P_1 \rightarrow P_2 
\rightarrow \cdots \rightarrow P_k \rightarrow P_{k+1} 
\rightarrow \cdots \] 
satisfying the following properties:
\begin{enumerate}
\item $G^{\mathrm{acy}}= 
\displaystyle\varinjlim P_k$, and 
$\iota_G \co  G \rightarrow G^{\mathrm{acy}}$ 
coincides with the limit map of the 
above sequence. 
\item $G \rightarrow P_k$ is a $2$--connected homomorphism. 
\end{enumerate}
\end{prop}
\noindent
From this proposition, we see, in particular, that the acyclic 
closure of a finitely presentable group is a countable set. 

\begin{remark}\label{rmk}
Here, we comment on the proofs of the above propositions. 
In the argument of the algebraic closure in \cite{le1}, 
Levine used the condition that a group $H$ is 
finitely normally generated by a subgroup $G$. In the case 
of the acyclic closure, we need the following 
alternative condition: the group $H$ is said to be 
{\it finitely homologically generated by a subgroup $G$} if 
\begin{enumerate}
\item The inclusion $G \to H$ induces a surjective homomorphism 
$H_1 (G) \to H_1 (H)$.
\item $H$ is generated by $G$ together with finite elements of $H$. 
\end{enumerate}
\noindent
As for the invisible subgroup, we need not change its definition. 
\end{remark}

\section{Structures of the groups $\Bn$ and 
$\Aut \Acy_n$}\label{ch:properties}
Using the results in the last section, 
we consider the acyclic closure 
$\iota_{F_n}\co F_n \to \Acy_n$ of $F_n$. 
Since the nilpotent completion $\Nil_n$ of $F_n$ is AC, 
there exists a unique homomorphism 
$p\co \Acy_n \to \Nil_n$ such that $p \circ \iota_{F_n}$ 
coincides with the natural map 
$F_n \to \Nil_n$, which is known to be injective. 
Hence $\iota_{F_n}$ is also injective. 

For each element $(G,\varphi_+,\varphi_-)$ of $\An$, 
we have a commutative diagram
\[\begin{CD}
F_n @>\varphi_->> G @<\varphi_+<< F_n\\ 
@V\iota_{F_n}VV @V\iota_{G}VV @VV\iota_{F_n}V \\
\Acy_n @>\cong>\varphi_-^{\mathrm{acy}}> G^{\mathrm{acy}} 
@<\cong<\varphi_+^{\mathrm{acy}}< \Acy_n 
\end{CD}\]
by \fullref{isom}. 
From this, we obtain a monoid homomorphism defined by
\[\Phi\co \An \longrightarrow \Aut \Acy_n \qquad 
\left( (G, \varphi_+, \varphi_-) \mapsto 
(\varphi_+^{\mathrm{acy}})^{-1} \circ 
\varphi_-^{\mathrm{acy}} \right)\]
and it induces a group homomorphism 
$\Phi\co \Bn \rightarrow \Aut \Acy_n$ by the commutativity of 
the diagram ($\ast$) in \fullref{ch:enlargement} 
whose homomorphisms are all 2--connected. 
\begin{thm}For each $n \ge 0$, the homomorphism 
$\Phi\co \Bn \to \Aut \Acy_n$ is an isomorphism.
\end{thm}
\begin{proof}
Assume that $(G, \varphi_+ , \varphi_-) \in \Ker \Phi$. Then 
$\varphi_+^{\mathrm{acy}}= \varphi_-^{\mathrm{acy}}:
\Acy_n \to G^{\mathrm{acy}}$, so that we have 
$\iota_G \circ \varphi_+ = 
\varphi_+^{\mathrm{acy}} \circ \iota_{F_n} = 
\varphi_-^{\mathrm{acy}} \circ \iota_{F_n} =
\iota_G \circ \varphi_-$. 
By \fullref{seq}, for large $k \ge 0$, we have 
$i_k \circ \varphi_+ = i_k \circ \varphi_- \co F_n \to P_k$ 
where $P_k$ is the $k^{\rm th}$ group of a sequence whose direct limit 
gives $G^{\mathrm{acy}}$, and $i_k \co G \to P_k$ is the composite 
of homomorphisms of the sequence from $G=P_0$ up to $P_k$. 
When we write 
$\varphi:=i_k \circ \varphi_+ = i_k \circ \varphi_-\co F_n \to P_k$, 
then $(G, \varphi_+ , \varphi_-) \in \An$ is homology cobordant 
to the identity element $(F_n, \Id, \Id)$ by a homology cobordism 
$G \stackrel{i_k}{\longrightarrow} P_k 
\xleftarrow{\varphi} F_n$. This shows 
that $\Phi$ is injective. 

On the other hand, given $\varphi \in \Aut \Acy_n$, we set 
$f:=\varphi \circ \iota_{F_n}\co F_n \to \Acy_n$. 
By \fullref{seq}, we have a 
sequence $\{ P_k \}$ of finitely presentable groups 
whose direct limit is $\Acy_n$. 
For large $k \ge 0$, we can take 
a lift $\widetilde{f}\co F_n \to P_k$ of $f$ with respect 
to the limit map $\iota\co P_k \to \Acy_n$, that is, we have 
$\iota \circ \widetilde{f} = f$. 
By definition, we have $\iota \circ i_k = \iota_{F_n}$ where 
$i_k\co F_n \to P_k$ is the composite of homomorphisms in the 
sequence. We can see that
$i_k$ and $\widetilde{f}$ are 2--connected homomorphisms, so that 
$(P_k, i_k, \widetilde{f})$ defines an element of $\An$. 
Taking their acyclic closures, we obtain 
$\Phi (P_k, i_k, \widetilde{f})=\varphi$. 
This completes the proof. 
\end{proof}

\begin{cor}\label{injective}
The monoid homomorphism $\Aut F_n \to \An$ 
and the group homomorphism $\Aut F_n \to \Bn \cong \Aut \Acy_n$ 
described in Section $\ref{ch:enlargement}$ 
are both injective.
\end{cor}
\begin{proof}
By the universal property of $\Acy_n$, two automorphisms of 
$\Acy_n$ are the same if and only if they coincide on 
the subgroup $F_n \subset \Acy_n$. 
The claim follows from this. 
\end{proof}

Hereafter we identify $\Bn$ with $\Aut \Acy_n$, and use 
only the latter. 
In the rest of this section, we describe 
some fundamental tools for understanding the structure 
of the group $\Aut \Acy_n$. 

\medskip
{\bf The Johnson homomorphisms}\qua
By Stallings' theorem, 
the inclusion $F_n \hookrightarrow \Acy_n$ induces 
isomorphisms on their nilpotent quotients. 
Therefore we have a natural homomorphism 
$\Phi_k\co \Aut \Acy_n \longrightarrow \Aut N_k$ 
for each $k \ge 2$. 
\begin{prop}\label{surjective}
For all $k \ge 2$, the homomorphisms 
$\Phi_k\co \Aut \Acy_n \to \Aut N_k$ are surjective. 
\end{prop}
\begin{proof}
Given an element $f \in \Aut N_k$, we denote 
by $\widetilde{f} \in \End F_n$ a lift 
of $f$. Since $\widetilde{f}$ induces an automorphism on $N_2$, 
$\widetilde{f}$ is a $2$--connected endomorphism. Then 
$\widetilde{f}^{\mathrm{acy}} \in \Aut \Acy_n$ is induced and it 
satisfies $\Phi_k (\widetilde{f}^{\mathrm{acy}})=f$.
\end{proof}
\noindent
By using $\Phi_k$, 
the Johnson homomorphism is defined as follows 
(see also \cite{mo} and \cite{ka}). We define a filtration of 
$\Aut \Acy_n$ by 
\[\Aut \Acy_n [1]:=\Aut \Acy_n , \qquad 
\Aut \Acy_n [k]:= \Ker \Phi_k \quad (k \ge 2). \]
On the other hand, we have an exact sequence
\[0 \longrightarrow \Hom (H_1 (F_n), (\Gamma^k F_n)/(\Gamma^{k+1} F_n) ) 
\longrightarrow \Aut N_{k+1} \longrightarrow \Aut N_k 
\longrightarrow 1,\]
where $(\Gamma^k F_n)/(\Gamma^{k+1} F_n)$ is known to be isomorphic 
to the degree $k$ part of the graded Lie algebra (over $\Z$) 
freely generated by 
the elements of $H_1 (F_n)$, so that $\Aut N_2$ acts on it. 
Explicitly, the isomorphism 
\begin{align*}
\Ker (\Aut N_{k+1} \to \Aut N_k) 
\we &
\Hom (H_1 (F_n), (\Gamma^k F_n)/(\Gamma^{k+1} F_n) )\\& 
=\Hom (F_n, (\Gamma^k F_n)/(\Gamma^{k+1} F_n) )
\end{align*}
is given by 
assigning to $f \in \Ker (\Aut N_{k+1} \to \Aut N_k)$ the 
homomorphism 
\[\begin{array}{ccc}
F_n \ni x_i & \mapsto & \widetilde{f}(x_i)x_i^{-1} \in 
(\Gamma^k F_n)/(\Gamma^{k+1} F_n)
\end{array}\]
where $\widetilde{f} \in \End F_n$ is a lift of $f$ and 
$\langle x_1,\ldots,x_n \rangle$ is 
a generating system of $F_n$. 
Note that this expression does not depend on the choices involved. 
If we define $J_k:=\Phi_{k+1} |_{\Aut \Acy_n [k]}$, 
we obtain an exact sequence
\[1 \to \Aut \Acy_n [k+1] \longrightarrow \Aut \Acy_n [k] 
\xrightarrow{J_k} 
\Hom (H_1 (F_n), (\Gamma^k F_n)/(\Gamma^{k+1} F_n) ) 
\to 1.\]
We call the homomorphism $J_k$ the $(k-1)^{\rm st}$ 
Johnson homomorphism. Note that $J_k$ is 
$\Aut \Acy_n$--equivariant, where $\Aut \Acy_n$ acts on 
$\Aut \Acy_n [k]$ by conjugation 
and acts on the target through $\Phi_2$.
\begin{example}
Consider the 2--connected endomorphism $\psi$ in \fullref{ex}. 
As an element of $\Aut \Acy_2$, $\psi$ belongs to $\Aut \Acy_2 [2]$. 
We calculate the image by the first Johnson 
homomorphism $J_2$. We write $H:=H_1 (F_2)$ and consider 
isomorphisms
\[\Hom (H_1 (F_2), (\Gamma^2 F_2)/(\Gamma^3 F_2)) \cong 
H^{\ast} \otimes (\Gamma^2 F_2)/(\Gamma^3 F_2) \cong 
H^{\ast} \otimes (\wedge^2 H).\]
Then we have 
$J_2 (\psi) = x_1^{\ast} \otimes (x_2 \wedge x_1)$.
\end{example}

\medskip
{\bf A refinement of the Johnson homomorphisms}\qua For each $k \ge 2$, 
we now give a refinement of the Johnson homomorphism 
whose target is abelian and bigger than that of the original. 
To construct the refinement, we need to fix a generating 
system $\langle x_1, \ldots, x_n \rangle$ of $F_n$. We show 
the following.
\begin{thm}\label{refinement}
For each $k \ge 2$, the Johnson homomorphism 
$J_k$ has a refinement 
\[\widetilde{J}_k \co  \Aut \Acy_n [k] \longrightarrow 
\Hom (H_1 (F_n), (\Gamma^k F_n)/(\Gamma^{2k-1} F_n) ) \]
whose target is also a finitely generated free abelian group. 
In fact, the composite with the natural projection 
$p_1\co (\Gamma^k F_n)/(\Gamma^{2k-1} F_n) \to 
(\Gamma^k F_n)/(\Gamma^{k+1} F_n)$ is the original Johnson 
homomorphism $J_k$. 
Moreover $\widetilde{J}_k$ is surjective, and 
the kernel of $\widetilde{J}_k$ coincides with 
$\Aut \Acy_n [2k-1]$, so that we have an exact sequence
\[1 \to \Aut \Acy_n [2k-1] \to \Aut \Acy_n [k] 
\xrightarrow{\widetilde{J}_k}
\Hom (H_1 (F_n), (\Gamma^k F_n)/(\Gamma^{2k-1} F_n) ) 
\to 1. \]
\end{thm}
\noindent
Note that for each $k \ge 2$, we have a direct sum decomposition 
\renewcommand{\theequation}{$\ast\ast$}
\begin{align}(\Gamma^k F_n)/(\Gamma^{2k-1} F_n) \cong 
\bigoplus_{j=k}^{2k-2}
\left( (\Gamma^j F_n)/(\Gamma^{j+1} F_n) \right) 
\end{align}
which is given by iterated extensions of 
$(\Gamma^k F_n)/(\Gamma^{k+1} F_n)$ by 
finitely generated\break free abelian groups 
$(\Gamma^i F_n)/(\Gamma^{i+1} F_n)$ for 
$k+1 \le i \le 2k-2$. 
Therefore\break $(\Gamma^k F_n)/(\Gamma^{2k-1} F_n)$ 
is also a finitely generated free abelian group. 
We also note that this direct sum decomposition is not canonical, 
except for the first projection 
$p_1\co  (\Gamma^k F_n)/(\Gamma^{2k-1} F_n) 
\to (\Gamma^k F_n)/(\Gamma^{k+1} F_n)$. 

The proof of \fullref{refinement} essentially uses 
the following. 
\begin{lem}\label{lem:refinement}
If we set $\Aut F_n[k]:=\Ker (\Aut F_n \to \Aut N_k)$, then 
\[\begin{array}{rccc}
\widetilde{J}_k\co & \Aut F_n[k] & \to &
\Hom (F_n, (\Gamma^k F_n)/(\Gamma^{2k-1} F_n) ) \\
& \mbox{\rotatebox[origin=c]{90}{$\in$}} & & 
\mbox{\rotatebox[origin=c]{90}{$\in$}} \\
& \varphi & \mapsto &
\left( \, x_i \mapsto \varphi (x_i) x_i^{-1} \, \right)
\end{array}\]
is a well-defined homomorphism.
\end{lem}
\begin{proof}
Given $\varphi$, $\psi \in \Aut F_n [k]$, we have
\[\widetilde{J}_k (\varphi \psi)(x_i) = 
\varphi(\psi(x_i)) x_i^{-1} = 
\varphi(\psi(x_i) x_i^{-1}) \cdot \varphi(x_i)x_i^{-1}.
\]
Since $\varphi(x_i) x_i^{-1}$, 
$\psi(x_i) x_i^{-1} \in \Gamma^k F_n$, 
and $(\Gamma^k F_n)/(\Gamma^{2k-1} F_n)$ is abelian, it suffices 
to show that $\Aut F_n [k]$ acts on 
$(\Gamma^k F_n)/(\Gamma^{2k-1} F_n)$ trivially. 
Every element $g \in (\Gamma^k F_n)/(\Gamma^{2k-1} F_n)$ can be 
written in a form $g= \prod_{i=1}^{l} 
[\cdots [g_{i1},g_{i2}],g_{i3}],\cdots],g_{ik}]$ where 
$g_{ij} \in F_n$, so that 
it suffices to show our claim in the case of\break
$g= [\cdots [g_{1},g_{2}],g_{3}],\cdots],g_{k}]$ where 
$g_{i} \in F_n$. 

Since $\varphi \in \Aut F_n [k]$, we see 
$\varphi(g)=[\cdots [g_1 r_1 ,g_2 r_2],g_3 r_3],\cdots],g_k r_k]$ 
for some $r_i \in \Gamma^k F_n$. 
We write $g^{(l)}:=[\cdots [g_{1},g_{2}],g_{3}],\cdots],g_{l}] 
\in \Gamma^l F_n$ for $2 \le l \le k$. Now we show 
that $g^{(l)} \equiv \varphi(g^{(l)}) \pmod{\Gamma^{k+l-1} F_n}$ 
by the induction on $l$. Our claim follows from it. 
For $l=2$, 
\begin{align*}
\varphi(g^{(2)}) &=[g_1 r_1, g_2 r_2]\\
&=[g_1 r_1, g_2] \cdot {}^{g_2}[g_1 r_1, r_2]\\
&={}^{g_1}[r_1, g_2] \cdot [g_1, g_2] \cdot {}^{g_2}[g_1 r_1, r_2]\\
& \equiv g^{(2)} \quad \pmod{\Gamma^{k+1} F_n}
\intertext{where we write ${}^a [b,c]$ for $a [b,c] a^{-1}$. 
When $g^{(i)} \equiv \varphi(g^{(i)}) \pmod{\Gamma^{k+i-1} F_n}$ 
follows for $2 \le i \le l$, we see}
\varphi(g^{(l+1)}) &=[\varphi(g^{(l)}), g_{l+1} r_{l+1}]\\
&=[\varphi(g^{(l)}),g_{l+1}] \cdot 
{}^{g_{l+1}}[\varphi(g^{(l)}), r_{l+1}]\\
&=[g^{(l)} r,g_{l+1}] \cdot 
{}^{g_{l+1}}[\varphi(g^{(l)}), r_{l+1}] \quad 
\mbox{for some $r \in \Gamma^{k+l-1}$}\\
&={}^{g^{(l)}}[r, g_{l+1}] \cdot 
[g^{(l)}, g_{l+1}] \cdot {}^{g_{l+1}}[\varphi(g^{(l)}),r_{l+1}]\\
& \equiv g^{(l+1)} \quad \pmod{\Gamma^{k+l} F_n}
\end{align*}
and this completes the proof.
\end{proof}
\noindent
By \fullref{lem:refinement}, we see that $\widetilde{J}_k$ 
gives a refinement of the 
Johnson homomorphism for $\Aut F_n$. 

\begin{proof}[Proof of \fullref{refinement}] 
If we restrict 
$\Phi_{2k-1}\co \Aut \Acy_n \to \Aut N_{2k-1}$ 
to the subgroup $\Aut \Acy_n [k]$, 
its image is contained in 
$\Ker(\Aut N_{2k-1} \to \Aut N_k)$. 
On the other hand, 
the map
\[\begin{array}{rccc}
\overline{J}_k\co & \Ker(\Aut N_{2k-1} \to \Aut N_k) & \to &
\Hom (F_n, (\Gamma^k F_n)/(\Gamma^{2k-1} F_n) ) \\
& \mbox{\rotatebox[origin=c]{90}{$\in$}} & & 
\mbox{\rotatebox[origin=c]{90}{$\in$}} \\
& f & \mapsto &
\left( \, x_i \mapsto \widetilde{f} (x_i) x_i^{-1} \, \right)
\end{array}
\]
where $\widetilde{f} \in \End F_n$ is a lift of $f$, defines 
a well-defined injective homomorphism 
by an argument similar to that 
in the proof of \fullref{lem:refinement}. 
Then we define a homomorphism $\widetilde{J}_k\co \Aut \Acy_n [k] 
\to \Hom (F_n, (\Gamma^k F_n)/(\Gamma^{2k-1} F_n) )$ by 
the composite
\[\Aut \Acy_n [k] \xrightarrow{\Phi_{2k-1}} 
\Ker(\Aut N_{2k-1} \to \Aut N_k) \xrightarrow{\overline{J}_k} 
\Hom (F_n, (\Gamma^k F_n)/(\Gamma^{2k-1} F_n) ).\]
It is easily checked that $\widetilde{J}_k$ gives 
a refinement of the Johnson homomorphism and that 
the kernel of $\widetilde{J}_k$ coincides with $\Aut \Acy_n [2k-1]$. 

To show that $\widetilde{J}_k$ is surjective, we recall 
the direct sum decomposition ($\ast\ast$). We write 
\[p_l\co(\Gamma^k F_n)/(\Gamma^{2k-1} F_n) \to 
(\Gamma^{k+l-1} F_n)/(\Gamma^{k+l} F_n),\quad (1 \le l \le k-1)\]
for the $l^{\rm th}$ projection. While each projection $p_l \; 
(2 \le l \le k-1)$ except $p_1$ 
is not given canonically, its restriction 
to $(\Gamma^{k+l-1} F_n)/(\Gamma^{2k-1} F_n)$ 
coincides with the natural projection 
$(\Gamma^{k+l-1} F_n)/(\Gamma^{2k-1} F_n) \to 
(\Gamma^{k+l-1} F_n)/(\Gamma^{k+l} F_n)$. 
Therefore, if we consider the isomorphism given by
\[\Hom(F_n, (\Gamma^k  F_n)/(\Gamma^{2k-1} F_n)) 
\xrightarrow[p_1 \oplus \cdots \oplus p_{k-1}]{\cong} 
\bigoplus_{j=k}^{2k-2}
\Hom \left( F_n,(\Gamma^j F_n)/(\Gamma^{j+1} F_n) \right),\]
the composite 
\[p_l \circ \widetilde{J}_k |_{\Aut \Acy_n [k+l-1]}\co
\Aut \Acy_n [k+l-1] \to 
\Hom ( F_n,(\Gamma^{k+l-1} F_n)/(\Gamma^{k+l} F_n) )\]
is nothing other than the original Johnson homomorphism $J_{k+l-1}$ 
for each $l=1,2,\ldots,k-1$. 
Since $J_k,\ldots,J_{2k-2}$ are all surjective, our claim follows. 
\end{proof}

\begin{remark}
The homomorphism $\widetilde{J}_k$ highly depends on the 
choice of a generating system of $F_n$, and $\widetilde{J}_k$ is 
not $\Aut \Acy_n$--equivariant for $k \ge 3$. 
This phenomenon is explained by 
using the Magnus expansion as follows. It is well known 
that the expansion of an element of $\Gamma^k F_n$ has a form 
of $1+(\mbox{degree $\ge k$--part})$. In terms of the Magnus 
expansion, our refinement $\widetilde{J}_k$ captures an 
information of the part from degree $k$ up to $(2k-2)$ of 
the expansion of $\widetilde{f} (x_i) x_i^{-1}$ 
under a fixed generating system of $F_n$. 
For a changing of a generating system, 
the Magnus expansion for each element intensively varies except 
that the first non-trivial homogeneous component in the positive 
degree part changes $\Aut \Acy_n$--equivariantly (see \cite{mo}, \cite{ka}). 
\end{remark}

{\bf The Magnus representation}\qua Here 
we define the Magnus representation for 
$\Aut \Acy_n$. While we call it the Magnus 
``representation'', it is actually a crossed homomorphism. 
The construction of the representation is based on 
Le Dimet's work \cite{ld}, 
where the Gassner representation 
for the pure braid group is extended to that for the 
string link concordance group. 

Before starting our discussion, we summarize our notation 
and rules. For a matrix $A$ with coefficients in a ring $R$, 
and a homomorphism 
$\varphi\co R \to R'$, we denote by ${}^{\varphi} A$ the 
matrix obtained from $A$ by applying $\varphi$ to each entry. 
When $R=\Z G$ for a group $G$ (or 
its Cohn localization mentioned below), 
we denote by $\overline{A}$ 
the matrix obtained from $A$ by applying the involution induced 
from $(x \mapsto x^{-1},\ x \in G)$ to each entry.

For a (finite) CW--complex $X$ and 
its regular covering $X_{\Gamma}$ with respect to a homomorphism 
$\pi_1 X \to \Gamma$, $\Gamma$ acts on $X_{\Gamma}$ from 
the right through its deck transformation group. Therefore 
we regard the $\Z \Gamma$--cellular chain complex 
$C_{\ast} (X_{\Gamma})$ of $X_{\Gamma}$ as a collection of 
free right $\Z \Gamma$--modules consisting of column vectors together 
with differentials given by left multiplications of matrices. 
For each $\Z \Gamma$--bimodule $A$, 
the twisted chain complex $C_{\ast} (X;A)$ is given 
by the tensor product of 
the right $\Z \Gamma$--module $C_{\ast} (X_{\Gamma})$ and 
the left $\Z \Gamma$--module $A$, 
so that $C_{\ast} (X;A)$ and 
$H_{\ast} (X;A)$ are right $\Z \Gamma$--modules.

To construct the Magnus representation for $\Aut \Acy_n$, we 
use the following special case of the {\it Cohn localization} 
(or the {\it universal localization}). 
We refer to Section 7 in \cite{co} for details.

\begin{prop}[Cohn \cite{co}]\label{prop:cohn} 
Let $G$ be a group and let $\varepsilon\co\Z G \to \Z$ be 
the augmentation map. Then there exists a pair of 
a ring $\Lambda_G$ and a ring homomorphism 
$l_G\co\Z G \to \Lambda_G$ satisfying the following properties:
\begin{enumerate}
\item For every matrix $m$ with coefficients in $\Z G$, 
if $\varepsilon(m)$ is invertible 
then $l_G (m)$ is also invertible.
\item The pair $(\Lambda_G, l_G)$ is universal among all pairs 
having the property 1. 
\end{enumerate}
Furthermore it is unique up to isomorphism. 
\end{prop}
\begin{example}
When $G=H_1 (F_n)$, we have 
\[\Lambda_G \cong \left\{ \, 
\displaystyle\frac{f}{g} \, \biggm| 
f,g \in \Z G, 
\varepsilon(g)=\pm 1 \right\}.\]
\end{example}

We write $x_i$ again for the image of $x_i$ by 
$\iota_{F_n}\co F_n=\langle x_1,x_2,\ldots,x_n \rangle 
\hookrightarrow \Acy_n$. 
\begin{prop}[Proposition 1.1 in \cite{ld}]\label{extfree} 
The homomorphism
\[\begin{array}{rccc}
v\co& \Lambda_{\Acy_n}^n & \to &
I(\Acy_n) \otimes_{\Acy_n} \Lambda_{\Acy_n} \\
& \mbox{\rotatebox[origin=c]{90}{$\in$}} & & 
\mbox{\rotatebox[origin=c]{90}{$\in$}} \\
& (a_1,\ldots,a_n) & \mapsto &
 \smash{\displaystyle\sum_{i=1}^n (x_i^{-1}-1) \otimes a_i}\vrule width 0pt depth15pt 
\end{array}\]
is an isomorphism of right $\Lambda_{\Acy_n}$--modules, where 
$I(\Acy_n):=\Ker (\varepsilon\co\Z \Acy_n \to \Z)$.
\end{prop}
\noindent
Note that each automorphism of $\Acy_n$ induces one of $\Z \Acy_n$. 
Moreover, by the universal property of $\Lambda_{\Acy_n}$, an 
automorphism of $\Lambda_{\Acy_n}$ is also induced. 

The proof of \fullref{extfree} is almost the same as that 
of Proposition 1.1 in \cite{ld}, once we show the following.
\begin{lem}
Let $G$ be a finitely presentable group, and let $f\co F_n \to G$ 
be a 2--connected homomorphism. Then 
$f_{\ast}\co H_i (F_n;\Lambda_G) \to H_i (G;\Lambda_G)$ is 
an isomorphism for $i =0,1,2$.
\end{lem}
\begin{proof}
We prove this lemma by using the idea of 
the proof of Proposition 2.1 in \cite{klw}. 
Let $X=K(F_n ,1)$ be a bouquet of $n$ circles 
and $Y=K(G,1)$ be a CW--complex 
constructed from a finite presentation of $G$. 
The number of cells of $Y$ up to degree 2 is finite. 
We denote by $f$ again for the continuous map 
from $X$ to $Y$ induced 
by the homomorphism $f\co F_n \to G$. 
Taking a mapping cylinder with respect to $f$, we obtain a 
CW--complex $M=K(G,1)$ where $X$ is contained as a subcomplex. 
The number of cells of $M$ up to degree 2 is also finite. 
Since $H_i (M,X)=0$ for 
$i=0,1,2$, we can take a partial chain homotopy 
$D_{i+1}\co C_i (M,X) \to C_{i+1} (M,X)$ of the partial chain complex 
$C_3 (M,X) \to \cdots \to C_0 (M,X) \to 0$ 
freely generated by relative cells of $(M,X)$. Namely, 
we have 
\begin{align*}
1 &= \partial_1 \circ D_1,\\
1 &= \partial_2 \circ D_2 + D_1 \circ \partial_1, \\
1 &= \partial_3 \circ D_3 + D_2 \circ \partial_2.
\end{align*}
Let $\widetilde{M}$ be the universal covering of $M$ and 
$\widetilde{X}$ be the inverse image of $X$ on $\widetilde{M}$. 
We choose a lift of each cell of $M$ on $\widetilde{M}$. 
Using the lifts of cells, we can define 
lifts $\widetilde{D}_{i+1}\co C_i (\widetilde{M},\widetilde{X}) 
\to C_{i+1} (\widetilde{M},\widetilde{X})$ of 
the chain homotopy $D_{i+1}$ for $i=0,1,2$, 
which are $\Z G$--equivariant. Then 
we define 
\begin{align*}
\Phi_0 &:= \widetilde{\partial_1} \circ \widetilde{D_1},\\
\Phi_1 &:= \widetilde{\partial_2} \circ \widetilde{D_2} + 
\widetilde{D_1} \circ \widetilde{\partial_1}, \\
\Phi_2 &:= \widetilde{\partial}_3 \circ \widetilde{D_3} + 
\widetilde{D_2} \circ \widetilde{\partial_2},
\end{align*}
\noindent
where $\partial_i$ are differentials of the chain complex 
$C_i (\widetilde{M},\widetilde{X})$. 
It is easily checked that 
$\Phi_i\co C_i(\widetilde{M},\widetilde{X}) \to 
C_i(\widetilde{M},\widetilde{X})$ $(i=0,1,2)$ gives 
a partial chain map, so that it 
induces a homomorphism 
$(\Phi_i)_{\ast}\co H_i(\widetilde{M},\widetilde{X}) \to 
H_i(\widetilde{M},\widetilde{X})$ for each $i=0,1,2$. 
Note that each $\Phi_i$ is a homomorphism between 
finitely generated free $\Z G$--modules which is the identity 
map on the base space. 
Then by the definition of the Cohn localization, 
\[\Phi_i \otimes_G 1\co C_i(\widetilde{M},\widetilde{X}) \otimes_G 
\Lambda_G \longrightarrow 
C_i(\widetilde{M},\widetilde{X}) \otimes_G 
\Lambda_G
\]
becomes an isomorphism for each $i=0,1,2$. 
Moreover $\Phi_i \otimes_G 1$ maps 
$\Ker (\widetilde{\partial_i} \otimes_G 1)$ 
onto itself, so that $(\Phi_i \otimes_G 1)$ induces an epimorphism 
on $H_i$.

On the other hand, since 
\begin{align*}
\Phi_0 \otimes_G 1 &= \widetilde{\partial_1} \circ \widetilde{D_1} 
\otimes_G 1,\\
\Phi_1 \otimes_G 1 &= (\widetilde{\partial_2} \circ \widetilde{D_2} + 
\widetilde{D_1} \circ \widetilde{\partial_1})
\otimes_G 1, \\
\Phi_2 \otimes_G 1 &= (\widetilde{\partial}_3 \circ \widetilde{D_3} + 
\widetilde{D_2} \circ \widetilde{\partial_2})
\otimes_G 1,
\end{align*}
\noindent
we see that $(\Phi_i \otimes_G 1)_{\ast}:
H_i(M,X;\Lambda_G) \to H_i(M,X;\Lambda_G)$ are 0--maps, and 
therefore $H_i (M,X;\Lambda_G)=0$ for $i=0,1,2$. Then 
$0=H_2 (M;\Lambda_G)=H_2 (G;\Lambda_G)$. From this, we see that
$f_{\ast}\co H_i (F_n;\Lambda_G) \to H_i (G,\Lambda_G)$ is an 
isomorphism for each $i=0,1,2$.
\end{proof}

\begin{definition}
For $1 \le i \le n$, we define a map 
$\partial / \partial x_i\co \Acy_n \to \Lambda_{\Acy_n}$ by 
\[\begin{array}{rccc}
\left( \displaystyle\frac{\partial}{\partial x_1}, 
\frac{\partial}{\partial x_2}, \ldots, 
\frac{\partial}{\partial x_n} \right):& \Acy_n & \to &
\Lambda_{\Acy_n}^n \\
& \mbox{\rotatebox[origin=c]{90}{$\in$}} & & 
\mbox{\rotatebox[origin=c]{90}{$\in$}} \\
& g & \mapsto & \fulloverline{v^{-1} ((g^{-1}-1) \otimes 1)}.
\end{array}\]
\end{definition}
\noindent
The above maps $\partial / \partial x_i$ coincide 
with the ordinary free differentials if we restrict them to $F_n$, and 
have similar properties. 
We refer to Proposition 1.3 in \cite{ld}. In particular, we have 
\[(g^{-1} -1) \otimes 1= \sum_{i=1}^n (x_i^{-1} -1) \otimes 
\fulloverline{\left(\frac{\partial g}{\partial x_i}\right)}\]
for any $g \in \Acy_n$ under our notation. 
\begin{definition}
We define the Magnus representation 
\[r\co \Aut \Acy_n \to M(n,\Lambda_{\Acy_n})\]
\[r(\varphi):=\left(\fulloverline{\left(\frac{\partial \varphi(x_j)}
{\partial x_i} \right)}\right)_{i,j}\leqno{\hbox{by setting}} \] for $\varphi \in \Aut \Acy_n$
\end{definition}
\begin{prop}
The Magnus representation $r$ is a crossed homomorphism. 
In particular, the image of $r$ is contained in the set 
of invertible matrices.
\end{prop}
\begin{proof}
For $\varphi$, $\psi \in \Aut \Acy_n$, we have
\begin{align*}
(\varphi\psi(x_j^{-1})-1) \otimes 1 &= 
\sum_{i=1}^n (x_i^{-1} -1) \otimes 
\fulloverline{\left(\frac{\partial \varphi\psi(x_j)}{\partial x_i}\right)}
\intertext{by definition. On the other hand,}
(\varphi\psi(x_j^{-1})-1) \otimes 1 
&= {}^{\varphi}((\psi(x_j^{-1})-1) \otimes 1) \\
&= \sideset{^{\varphi}\!}{}
{\Tmatrix{\sum_{k=1}^{n} (x_k^{-1} -1) \otimes 
\fulloverline{\left(\frac{\partial \psi (x_j)}{\partial x_k}\right)}}}\\
&= \sum_{k=1}^{n} (\varphi(x_k^{-1})-1) \otimes 
\sideset{^{\varphi\;}\!}{}
{\mathop{\fulloverline{
\left( \frac{\partial \psi(x_j)}{\partial x_k} \right) 
}}\nolimits} \\
&= \sum_{k=1}^{n} 
\left\{ (\varphi(x_k^{-1})-1) \otimes 1 \right\} 
\cdot \sideset{^{\varphi\;}\!}{}
{\mathop{\fulloverline{
\left( \frac{\partial \psi(x_j)}{\partial x_k} \right) 
}}\nolimits} \\
&= \sum_{k=1}^{n} 
\left\{ \sum_{i=1}^{n} (x_i^{-1} -1) \otimes 
\fulloverline{\left(\frac{\partial \varphi(x_k)}{\partial x_i}\right)}
\right\} \cdot \sideset{^{\varphi\;}\!}{}
{\mathop{\fulloverline{
\left( \frac{\partial \psi(x_j)}{\partial x_k} \right) 
}}\nolimits} \\
&= \sum_{i=1}^{n} (x_i^{-1} -1) \otimes \left\{\sum_{k=1}^n 
\fulloverline{\left( \frac{\partial \varphi(x_k)}{\partial x_i}
\right) } \cdot 
\sideset{^{\varphi\;}\!}{}
{\mathop{\fulloverline{
\left( \frac{\partial \psi(x_j)}{\partial x_k} \right) 
}}\nolimits}\right\}.
\intertext{Hence we obtain}
\fulloverline{\frac{\partial \varphi\psi(x_j)}{\partial x_i}}
&= \sum_{k=1}^n 
\fulloverline{\left( \frac{\partial \varphi(x_k)}{\partial x_i}
\right) } \cdot 
\sideset{^{\varphi\;}\!}{}
{\mathop{\fulloverline{
\left( \frac{\partial \psi(x_j)}{\partial x_k} \right) 
}}\nolimits}
\end{align*}
\noindent
which shows that 
$r(\varphi \psi) = r(\varphi) \cdot {}^{\varphi} r(\psi)$.
\end{proof}
\noindent
Note that the composite $\Z F_n \xrightarrow{\iota_{F_n}} \Z \Acy_n 
\xrightarrow{l_{\Acy_n}} \Lambda_{\Acy_n}$ 
is injective, for the composite of the ring homomorphism
$\Z \Acy_n \to \Z \Nil_n$ with the Magnus expansion, 
which can be extended to $\Z \Nil_n$ and is injective on $\Z F_n$, 
satisfies the property 1 of \fullref{prop:cohn}, so that 
the Magnus expansion is extended for $\Lambda_{\Acy_n}$. 
Hence the Magnus representation defined 
here certainly gives a generalization of the original 
$r \co \Aut F_n \to GL(n,\Z F_n)$. 
\begin{example}\label{twoconn}
Consider the 2--connected endomorphism $\psi$ in \fullref{ex}. 
Then 
\[r(\psi)=
\left(\begin{array}{cc}
1+x_2^{-1}x_1^{-1}-x_1 x_2 x_1^{-1}x_2^{-1} x_1^{-1} & 0 \\
x_1^{-1}-x_2 x_1^{-1}x_2^{-1}x_1^{-1} & 1
\end{array}\right).\]
Reducing the coefficients to 
$\Lambda_{H_1 (\Acy_2)}=\Lambda_{H_1 (F_2)}$, we obtain the matrix
\[\left(\begin{array}{cc}
1+x_1^{-1}x_2^{-1}-x_1^{-1} & 0 \\
x_1^{-1}-x_1^{-2} & 1
\end{array}\right)\]
whose determinant is $1+x_1^{-1}x_2^{-1}-x_1^{-1}$. 
Since this value is not invertible in $\Z H_1 (F_2)$, 
we see that $\psi \notin \Aut F_2$.
\end{example}

\section{Main results}\label{ch:results}
In this section, we return to our discussion on 
homology cylinders, and consider some relationships between 
homology cylinders and the acyclic closure of a free group.

For each homology cylinder $(M,i_+,i_-) \in \Cg$, 
we obtain a commutative diagram
\[\begin{CD}
F_{2g} @>i_->> \pi_1 M @<i_+<< F_{2g}\\ 
@V\iota_{F_{2g}}VV @V\iota_{\pi_1 M}VV @VV\iota_{F_{2g}}V \\
\Acy_{2g} @>\cong>i_-^{\mathrm{acy}}> (\pi_1 M)^{\mathrm{acy}} 
@<\cong<i_+^{\mathrm{acy}}< \Acy_{2g} 
\end{CD}\]
by \fullref{isom}. 
From this, we obtain a monoid homomorphism defined by
\[\sigma^{\mathrm{acy}}\co \Cg \longrightarrow \Aut \Acy_{2g} \qquad 
\left( \, (M,i_+,i_-) \mapsto 
(i_+^{\mathrm{acy}})^{-1} \circ i_-^{\mathrm{acy}} \, \right)\]
and it induces a group homomorphism 
$\sigma^{\mathrm{acy}}\co \Hg \rightarrow \Aut \Acy_{2g}$. 

Our first result is 
a generalization of Dehn--Nielsen's theorem. 
Recall that $\zeta \in F_{2g} \subset \Acy_{2g}$ 
is a word corresponding to the boundary loop of $\Sg$. 
\begin{thm}\label{autsp}
The image of $\sigma^{\mathrm{acy}}\co \Hg \to \Aut \Acy_{2g}$ is 
\[\Aut_0 \Acy_{2g} :=
\{ \varphi \in \Aut \Acy_{2g} \mid
\varphi(\zeta) = \zeta \in \Acy_{2g} \}.\]
\end{thm}

\begin{proof} 
By the definition of homology cylinders, we have 
$i_+ (\zeta)=i_- (\zeta) \in \pi_1 M$ for every homology cylinder 
$(M,i_+,i_-)$. Hence we see that the image of 
$\sigma^{\mathrm{acy}}$ 
is contained in $\Aut_0 \Acy_{2g}$. 

Conversely, given an element $\varphi \in 
\Aut_0 \Acy_{2g}$, 
we construct a homology cylinder $M=(M,i_+,i_-)$ satisfying 
$\sigma^{\mathrm{acy}}(M)=\varphi$. 
The construction is based on that of 
\fullref{thm:gl} (which is 
Theorem 3 in \cite{gl}). 
In our context, however, we must pay an extra attention because 
we do not have a result which directly corresponds to 
Lemma 4.6 in \cite{gl}. 
This occurs from the fact that 
$\iota_{F_{2g}}\co F_{2g} \to \Acy_{2g}$ is injective although 
the composite 
$F_{2g} \to \Acy_{2g} \to \Acy_{2g}/(\Gamma^k \Acy_{2g}) 
\cong N_k$ is surjective. 

We begin by taking two continuous maps 
$f_+, f_-\co \Sg \to K(\Acy_{2g},1)$ corresponding to 
homomorphisms $\iota_{F_{2g}}, \varphi \circ 
\iota_{F_{2g}}\co F_{2g} \to \Acy_{2g}$, 
respectively. Since $\iota_{F_{2g}}(\zeta)=
\varphi \circ \iota_{F_{2g}}(\zeta)$, we can combine the 
two maps and obtain a map $f:=f_+ \cup f_- \co  
\Sigma_{2g}=\Sg \cup (-\Sg) \to K(\Acy_{2g},1)$. The pair 
$(\Sigma_{2g}, f)$ defines an element of the second bordism 
group $\Omega_2 (\Acy_{2g})$ of $K(\Acy_{2g},1)$, which is 
naturally isomorphic to $H_2 (\Acy_{2g})$. Since 
$H_2 (\Acy_{2g})=0$ as mentioned in \fullref{2cn}, 
there exist a compact oriented 
3--manifold $M$ whose boundary is $\Sigma_{2g}$, 
and a map $\Phi\co M \to K(\Acy_{2g},1)$ extending the map $f$. We 
write $i_+,i_-\co \Sg \to \partial M$ for embeddings onto 
domains of $f_+,f_-$, respectively. 

Since $\iota_{F_{2g}}$ is 2--connected, 
we have $H_1 (M) \cong i_+(H_1 (\Sg)) \oplus \Ker \Phi_\ast$. 
If $\Ker \Phi_\ast=0$, then $i_+,i_-\co H_1 (\Sg) \to H_1 (M)$ are 
both isomorphisms. In particular, $H_1 (M,\partial M)=0$ and 
therefore $H^1 (M,\partial M) \cong H_2 (M) =0$, so that 
the triplet $(M,i_+,i_-)$ gives a homology cylinder satisfying 
$\sigma^{\mathrm{acy}}(M)=\varphi$. 
If not, we perform surgery on the map $\Phi$ to kill 
$\Ker \Phi_\ast$. 

Now we take a non-trivial element $\alpha \in \Ker \Phi_\ast$. 
If there exists a representative $C \in \pi_1 M$ 
by a simple closed curve such that $\Phi (C)=1 \in \Acy_{2g}$, 
then we can do surgery on $C$ and extend $\Phi$ over the trace 
of the surgery. However, we cannot necessarily take such a 
simple closed curve. Then we replace $(M,i_+,i_-)$ by a 
manifold which is homology bordant to $M$ over $K(\Acy_{2g},1)$ 
and for which we can take a simple closed curve 
which represents $\alpha$ and whose 
image by $\Phi$ is trivial in $\Acy_{2g}$. A construction 
of such a homology bordant manifold is given as follows. 

For the induced homomorphism $i_+\co F_{2g} \to \pi_1 M$, 
by the universal property of acyclic closures, we have 
a homomorphism $i_+^{\mathrm{acy}}\co \Acy_{2g} \to 
(\pi_1 M)^{\mathrm{acy}}$ satisfying $i_+^{\mathrm{acy}} \circ 
\iota_{F_{2g}}=\iota_{\pi_1 M} \circ i_+$. Similarly we 
have $\Phi^{\mathrm{acy}}\co (\pi_1 M)^{\mathrm{acy}} \to \Acy_{2g}$ 
satisfying $\Phi=\Phi^{\mathrm{acy}} \circ \iota_{\pi_1 M}$. Then 
$\Phi^{\mathrm{acy}} \circ i_+^{\mathrm{acy}} \circ \iota_{F_{2g}} 
= \Phi^{\mathrm{acy}} \circ \iota_{\pi_1 M} \circ i_+ = 
\iota_{F_{2g}}$, 
so that $\Phi^{\mathrm{acy}} \circ i_+^{\mathrm{acy}} 
= \Id_{\Acy_{2g}}$. 
In particular, $\Phi^{\mathrm{acy}}$ is onto.

Take a simple closed curve $C$ representing $\alpha \in 
\Ker \Phi_{\ast}$. Since $\Phi_\ast (\alpha)=0 \in H_1 (\Acy_{2g})$, 
we can write $\Phi (C)=\prod_{i=1}^{l} [h_{i1},h_{i2}]$ 
with $h_{ij} \in \Acy_{2g}$. We take an acyclic system 
\[S \, : \, x_i=w_i (x_1,x_2,\ldots,x_m) \quad (i=1,2,\ldots,m)\]
over $\pi_1 M$ whose solution in $(\pi_1 M)^{\mathrm{acy}}$ 
contains 
\[\{ i_+^{\mathrm{acy}}(h_{11}), i_+^{\mathrm{acy}}(h_{12}),\ldots, 
i_+^{\mathrm{acy}}(h_{l1}), i_+^{\mathrm{acy}}(h_{l2}) \}.\] 
We attach a 1--handle to $M \times \{ 1 \} \subset M \times I$ 
for each variable $x_i$ and write $x_i$ again for the added 
generator on the fundamental group of the resulting cobordism. 
We also attach a 2--handle along 
the loop $x_i w_i^{-1}$ for each $i=1,2,\ldots,m$ with any framing. 
We denote the resulting cobordism by $X_S$. Then 
\[\pi_1 X_S=(\pi_1 M \ast \langle x_1,\ldots,x_m \rangle) 
\big/ (x_1 w_1^{-1},\ldots,x_m w_m^{-1}).\]
Let $M_S$ be another part of the boundary of $X_S$, namely 
$M \cup M_S =\partial X_S$ and 
$M \cap M_S =\partial M = \partial M_S$. We define a homomorphism 
$\Phi_S\co  \pi_1 X_S \to (\pi_1 M)^{\mathrm{acy}}$ which lifts 
$\iota_{\pi_1 M}$ by sending $x_i$ to the corresponding solution 
of $S$. Then the composite $\Phi^{\mathrm{acy}} \circ \Phi_S\co 
\pi_1 X_S \to \Acy_{2g}$ 
gives a continuous map $\Phi_S \co X_S \to K(\Acy_{2g},1)$ which 
extends $\Phi\co M \to K(\Acy_{2g},1)$. 

Since $\iota_{\pi_1 M}$ is 2--connected and has a lift $\Phi_S$, 
the homomorphism $H_1 (M) \to H_1 (X_S)$ induced from 
the inclusion $M \hookrightarrow X_S$ is injective. On the other hand, 
by the definition of $X_S$, this homomorphism is onto, 
hence an isomorphism. We can also see that $H_2 (M) \to H_2 (X_G)$ is 
onto, for, in terms of chain complexes consisting of handles, 
the boundary of each newly added 2--handle corresponding to the 
relation $x_i w_i^{-1}$, has a form of $[x_i]+
(\mbox{1--handles in $M \times I$})$ which shows that we have 
a surjective homomorphism from the module of 2--cycles on $M \times I$ 
to that on $X_S$. 

We also show that $H_1 (M_S) \to H_1 (X_S)$ is an isomorphism. 
The surjectivity follows by considering the dual handle 
decomposition of $X_S$, namely $X_S$ is constructed 
from $M_S \times I$ by attaching 
3-- and 2--handles. To show the injectivity, we see 
that $H_2 (X_S,M_S)=0$. By the Poincar\'e--Lefschetz duality, 
$H_2 (X_S,M_S) \cong H^2 (X_S,M)$. On the other hand, we have 
an exact sequence 
\[\begin{CD}
H^1 (X_S) @>>> H^1 (M) @>>> H^2 (X_S,M) @>>> H^2 (X_S) @>>> H^2 (M). 
\end{CD}\]
Since $H_\ast (M) \to H_\ast (X_S)$ is an isomorphism for 
$\ast =1$ and onto for $\ast=2$, the first map 
is an isomorphism, and the last one is injective. Hence 
$H^2 (X_S,M)=0$. 

From the above argument, we see that $M$ and $M_S$ are homology 
bordant over $K(\Acy_{2g},1)$ by 
the bordism $X_S$ and $\Phi_S$. Note that this bordism preserves 
the direct sum decomposition 
$H_1 (M) \cong i_+(H_1 (\Sg)) \oplus \Ker \Phi_\ast$, namely 
we also have $H_1 (M_S) \cong i_+(H_1 (\Sg)) \oplus \Ker \Phi_{G \ast}$ 
and we can take $\overline{\alpha} \in \Ker \Phi_{G \ast}$ which 
corresponds to $\alpha$. 
Consider the simple closed curve $C$ taken at the beginning of 
this argument. Since $\pi_1 M_S \to \pi_1 X_S$ is onto, there 
exists a simple closed curve $\overline{C}$ which attains 
$C$ in $\pi_1 X_S$. Now $h_{ij} \in \Acy_{2g}$ are contained 
in the image of $\Phi_S\co \pi_1 M_S \to \Acy_{2g}$, so that we 
can take $\wwbar{h_{ij}} \in \pi_1 M_S$ attaining $h_{ij}$. 
Then the simple closed curve $\overline{C} \left( 
\prod_{i=1}^l [\wwbar{h_{i1}},\wwbar{h_{i2}}] \right)^{-1}$ 
represents $\overline{\alpha}$ and is mapped by $\Phi_S$ 
to the trivial element of $\Acy_{2g}$. 

The rest of the proof is almost the same as Theorem 3 in \cite{gl}. 
The difference is that our killing $\Ker \Phi_{\ast}$ goes, 
if necessary, with changing the manifold to 
some homology bordant one at each step. 
\end{proof}

Next, we consider $\Hg[[\infty]]:= \Ker 
\left( \sigma^{\mathrm{acy}}\co  \Hg \to \Aut_0 \Acy_{2g} \right)$. 
In contrast to the case of the mapping class group, 
$\Hg[[\infty]]$ is non-trivial. 
Indeed the homology cobordism group $\Theta_{\Z}^3$ of homology 
3--spheres is contained in it. Our second result gives an 
additive invariant for $\Hg[[\infty]]$. 

For any element $[(M,i_+,i_-)] \in \Hg[[\infty]]$, we have 
$i_+^{\mathrm{acy}}=i_-^{\mathrm{acy}}\co 
\Acy_{2g} \to (\pi_1 M)^{\mathrm{acy}}$. Consider 
the composite $f\co M \to K(\Acy_{2g},1)$ of continuous maps 
\[\begin{CD}
M @>>> K(\pi_1 M,1) @>>> K((\pi_1 M)^{\mathrm{acy}},1) 
@>>> K(\Acy_{2g},1)
\end{CD}\]
where the last map is induced by 
the homomorphism $(i_+^{\mathrm{acy}})^{-1}=
(i_-^{\mathrm{acy}})^{-1}$. After adjusting by some homotopy, 
if necessary, we have 
$f \circ i_+ = f \circ i_-\co \Sg \to K(\Acy_{2g},1)$. 
Then we can consider a closed 3--manifold 
\[C_M =M/(i_+ (x)= i_-(x)) \quad x \in \Sg\]
and a continuous map $\widehat{f}\co C_M \to 
K(\Acy_{2g},1)$ induced from $f$. 
Note that $C_M$ is also obtained by gluing 
$1_{\Cg}=\Sg \times I$ along their boundaries. 
We call this operation {\it closing}. 
\begin{thm}\label{invariant}
The map $\theta\co \Hg[[\infty]] \to H_3 (\Acy_{2g})$ given by
\[\theta ([M,i_+,i_-]) := 
\widehat{f} ([C_M ]) \quad \in H_3 (\Acy_{2g})\]
is a well-defined homomorphism. Moreover it is surjective.
\end{thm}
\begin{proof}
The proof is divided into three steps:
\begin{itemize}
\item $\widehat{f} ([C_M])$ depends only on 
the homology cobordism class of $(M,i_+,i_-)$, 
\item $\theta$ is actually a homomorphism,
\item $\theta$ is onto. 
\end{itemize}

{\bf Step 1}\qua Let $(M,i_+,i_-)$ and $(N,j_+,j_-)$ be 
homology cylinders contained in the same homology 
cobordism class in $\Hg[[\infty]]$. 
We write $f_M\co M \to K(\Acy_{2g},1)$ and 
$f_N\co N \to K(\Acy_{2g},1)$ for maps constructed 
as above, respectively. 
We take a homology cobordism 
$W$ satisfying 
\[\partial W = M \cup (-N)/(i_+ (x)= j_+(x) , \,
i_- (x)=j_-(x)) \quad x \in \Sg,\]
so that we have $i_+=j_+,i_-=j_-\co \Sg \to W$ and 
\[i_+^{\mathrm{acy}}=j_+^{\mathrm{acy}}=
i_-^{\mathrm{acy}}=j_-^{\mathrm{acy}}\co \Acy_{2g} \to 
(\pi_1 W)^{\mathrm{acy}}.\]
Note that, by Stallings' theorem,
\[\Acy_{2g} \cong (\pi_1 M)^{\mathrm{acy}} \cong 
(\pi_1 W)^{\mathrm{acy}} \cong (\pi_1 N)^{\mathrm{acy}}.\]
Then the homomorphism $(i_+^{\mathrm{acy}})^{-1} \circ 
\iota_{\pi_1 W} (= (j_+^{\mathrm{acy}})^{-1} \circ 
\iota_{\pi_1 W})$ induces a map $f_W\co W \to K(\Acy_{2g},1)$ 
extending $f_M$ and $f_N$. 

We take a closed tubular neighborhood $V$ of 
$\partial M=\partial N$ in $\partial W$. $V$ is 
diffeomorphic to $\Sigma_{2g} \times I$. We glue 
$1_{\Cg} \times I=(\Sg \times I) \times I$ to $W$ by 
identifying $(\partial 1_{\Cg}) \times I$ with $V$. 
The resulting 4--manifold $\widehat{W}$ gives a 
homology cobordism between $C_M$ and $C_N$. 
Moreover, for $f_M \circ i_+ = f_N \circ j_+$ and 
$f_M \circ i_- = f_N \circ j_-$ are homotopic, the map 
$f_W$ is extended to a map 
$\widehat{f}_{W}\co \widehat{W} \to K(\Acy_{2g},1)$ whose 
restriction to $C_M$ (resp.~$C_N)$ coincides with  
$\widehat{f}_M$ (resp.~$\widehat{f}_N$).
Therefore we see 
$\widehat{f}_M ([C_M])=\widehat{f}_N ([C_N])$.

\medskip
{\bf Step 2}\qua For $[(M,i_+,i_-)]$, $[(N,j_+,j_-)] \in 
\Hg[[\infty]]$, we construct a 4--manifold 
\[W=(M \times I) \cup (N \times I) \cup 
(1_{\Cg} \times [0,3])\]
by the following gluing rule. We decompose 
$\partial (1_{\Cg} \times [0,3])$ into 
\[((\partial 1_{\Cg}) \times [0,1] )\cup
((\partial 1_{\Cg}) \times [1,2] )\cup
((\partial 1_{\Cg}) \times [2,3] )\cup
(1_{\Cg} \times \{ 0,3 \}).\]
We glue $(\partial 1_{\Cg}) \times [0,1]$ to 
$(\partial M) \times I$. We also glue 
$(\partial 1_{\Cg}) \times [2,3]$ to 
$\partial N \times I$ with opposite direction of unit intervals 
and opposite markings of homology cylinders $N$ and $1_{\Cg}$. 
Namely, for example, 
we identify $(\Sg \times \{1\}) \times \{0\}$ with 
$i_+ (\Sg) \times \{ 0 \} \subset M \times I$, and 
$(\Sg \times \{0\}) \times \{3\}$ with 
$j_+ (\Sg) \times \{ 0 \} \subset N \times I$. 
Then $\partial W$ consists of $C_M$, $C_N$, 
and 
\[(-M \times \{1 \}) \cup (-N \times \{1 \}) \cup 
((\partial 1_{\Cg}) \times [1,2]).\]
The last component is nothing other than $-C_{M \cdot N}$. 
Hence $W$ is a cobordism between $C_{M} \sqcup C_{N}$ 
and $-C_{M \cdot N}$. 

Since $f_M \circ i_+$, $f_M \circ i_-$, 
$f_N \circ j_+$ and $f_N \circ j_-$ are all homotopic, 
$(f_M \times I)\cup(f_N \times I)\co 
(M \times I) \cup (N \times I) \to K(\Acy_{2g},1)$ is 
extended to a map $f_W\co W \to K(\Acy_{2g},1)$ whose 
restriction to $C_{M}$, $C_{N}$ and 
$-C_{M\cdot N}$ coincide with  
$\widehat{f}_M$, $\widehat{f}_N$ and $\widehat{f}_{M \cdot N}$, 
respectively. Hence we see that 
$\theta (M) + \theta (N) = \theta(M \cdot N)$. 

\medskip
{\bf Step 3}\qua We construct a homology cylinder which attains 
a given element $\alpha \in H_3 (\Acy_{2g})$ by $\theta$. 
First, we use a construction used in \cite{orr}, \cite{le1}. 

Take a bouquet of $2g$ circles tied at the base point $\ast$ as 
a model of $K(F_{2g},1)$, and take a CW--complex realizing 
$K(\Acy_{2g},1)$. 
The homomorphism $\iota_{F_{2g}}\co F_{2g}=
\langle \gamma_1, \ldots, \gamma_{2g} \rangle 
\to \Acy_{2g}$ induces 
a continuous map $I_{F_{2g}}\co K(F_{2g},1) \to K(\Acy_{2g},1)$. 
As usual, we can assume that 
$K(F_{2g},1)$ is a subcomplex of $K(\Acy_{2g},1)$ having only one 
0--cell $\ast$, namely $I_{F_{2g}}$ is an inclusion. 
Let $\widehat{K}$ be the mapping cone of $I_{F_{2g}}$, which 
is obtained from the CW--complex $K(\Acy_{2g},1)$ 
by attaching a 2--cell $e_i$ along each loop representing 
$\gamma_i \in F_{2g} \subset \Acy_{2g}$ for $i=1,\ldots,2g$. 
We take an interior point $p_i$ in each 2--cell $e_i$. 

Since $H_3 (\Acy_{2g}) = H_3(K(\Acy_{2g},1)) \cong H_3 (\widehat{K}) \cong 
\Omega_3 (\widehat{K})$, we can take a pair of 
a closed oriented 3--manifold $M$ and 
a continuous map $f\co M \to \widehat{K}$ corresponding to $\alpha$ 
by this isomorphism. By applying a method used in the proof of 
Theorem 4 in \cite{orr} or Theorem 2 in \cite{le1}, 
we can construct a pair of 
a closed oriented 3--manifold $M'$ and 
a continuous map $f'\co M' \to \widehat{K}$ satisfying 
the following.
\begin{itemize}
\item $(M',f')$ is bordant to $(M,f)$ over $\widehat{K}$, that is, 
$f'([M'])=f([M]) \in H_3 (\widehat{K})$.
\item $p_i$ is a regular value of $f'$ for each $i=1,\ldots,2g$.
\item $L_i:=(f')^{-1}(p_i)$ is not empty and is 
a knot in $M'$ for each $i=1,\ldots,2g$.
\item There exists a closed tubular neighborhood $T_i$ of $L_i$ 
such that the diagram 
\[\begin{CD}
T_i @>f'>> \widehat{K} \\
@V{h_i}V{\cong}V @AAg_iA \\
S^1 \times D^2 @>\mathrm{proj}>> D^2
\end{CD}\]
commutes where $h_i$ is a homeomorphism sending $L_i \subset 
T_i$ onto $S^1 \times \{ 0 \}$, and $g_i$ is a homeomorphism 
onto the 2--cell $e_i$ sending $1 \in D^2$ to the base point $\ast$
of $\widehat{K}$.
\end{itemize}
\noindent
For simplicity, we write $(M,f)$ again for $(M',f')$. 
We take a meridian loop 
$m_i:=h_i^{-1} (\{1 \} \times \partial D^2)$ and 
a point $q_i:= h_i^{-1} (1,1)$ 
on $m_i$ for each $i$. Note that $f(q_i)=\ast$. 
We orient $m_i$ by using the orientation of the loop 
representing $\gamma_i$.

We take a base point $q \in M \setminus (\bigcup_{i=1}^{2g} T_i) 
\subset M$. By using some homotopy, we can assume $f(q)=\ast$. 
We connect the points $q_{i}$ and $q_{g+i}$ by a path $l_i^1$ in 
$M \setminus (\bigcup_{i=1}^{2g} T_i)$ such that $f$ is the constant 
map to $\ast$ on a neighborhood of $l_i^1$. We also connect 
a point of $l_i^1$ and $q$ by a path $l_i^2$ satisfying 
the same condition for $f$ as $l_i^1$. The condition for $f$ on the 
neighborhood of $l_i^1$, $l_i^2$ will be satisfied 
by surgeries using 1--handles. The meridian loops $m_i$ and 
the paths $l_i^1$, $l_i^2$ form 
a graph with trivalent vertices except one $g$--valent 
vertex $q$ as depicted below.

$$\includegraphics{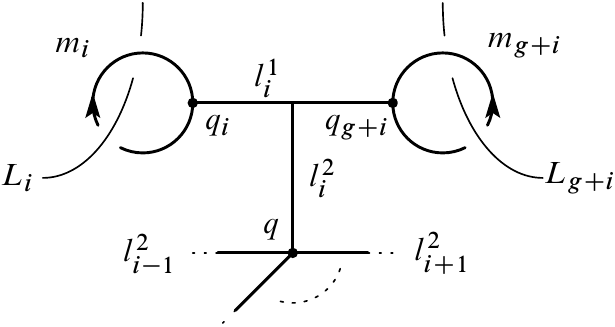}$$
In a tubular neighborhood of this graph, we take an embedded ribbon 
whose boundary is connected and is sent by $f$ to the loop 
$\zeta=\prod_{i=1}^g [\gamma_i,\gamma_{g+i}]$ in $\widehat{K}$ as follows. 
$$\includegraphics{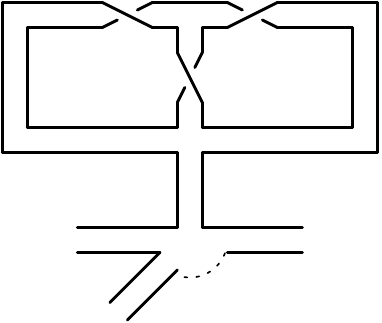}$$
We fatten the ribbon to construct an embedding of 
$1_{\Cg}=\Sg \times I$, in which the boundary of the ribbon divides 
$\partial 1_{\Cg} \cong \Sg \cup (-\Sg)$ into two embeddings 
$i_+$, $i_-$ of $\Sg$. 

We perform surgery at $L_1, L_2, \ldots, L_{2g}$ with framings 
given by $(h_i)^{-1} (S^1 \times \{1 \})$, on which $f$ is the 
constant map to $\ast$, 
and denote the resulting manifold by $M_L$. By the choice of 
framings, there exists a map $f_L\co M_L \to K(\Acy_{2g},1) 
\subset \widehat{K}$ which is bordant to $(M,f)$ over 
$\widehat{K}$. Note that 
$f_L ([M_L]) = \alpha \in H_3 (K(\Acy_{2g},1))$. 

We remove the embedded $1_{\Cg}$ from $M_L$ to construct 
a 3--manifold $N$ with boundary. $\partial N$ is diffeomorphic 
to $\Sg \cup (-\Sg)$, and has an orientation preserving embedding 
$i_+$ and a reversing one $i_-$. Now 
$f_L \circ i_+ = f_L \circ i_- = \iota_{F_{2g}}\co  \pi_1 \Sg \cong 
F_{2g} \to \Acy_{2g}$, so that $f_L \circ i_+$ and 
$f_L \circ i_-$ are 2--connected homomorphisms preserving 
$\zeta$. Only we have to do is to construct a homology cylinder 
from the data $N$, $i_+,i_-\co \Sg \to \partial N$ 
and $f_L\co N \to K(\Acy_{2g},1)$ with 
keeping its bordism class over $K(\Acy_{2g},1)$, which is 
the same situation of the proof of \fullref{autsp}. 
Using the same argument, we obtain a desired homology cylinder.
\end{proof}
\begin{remark}
As seen in Step 1 of the above proof, two 
homology cylinders belonging to the same class in $\Hg$ give 
homology cobordant closed 3--manifolds by closing.
\end{remark}
\begin{remark}
We do not know, at present, whether $H_3 (\Acy_{2g})$ is trivial or not. 
This situation is similar to 
that of some (string) link concordance invariant using the algebraic 
closure of a free group in \cite{le1} and \cite{le2}. 
It is easily checked that the homomorphism $\theta$ 
is trivial if we restrict it to 
$\Theta_{\Z}^3 \subset \Hg[[\infty]]$. 

We also have the following important problem. 
It is not known that $\Hg[[\infty]]$ 
coincides with the group 
\begin{align*}
\Hg[\infty]:=&\bigcap_{k}\Ker 
\left( \sigma_k\co  \Hg \to \Aut N_{k} \right)\\
=&\Ker 
\left( \sigma^{\mathrm{nil}}\co  \Hg \to \Aut \Nil_{2g} \right),
\end{align*}
\noindent
in which $\Hg[[\infty]]$ is contained. 
This problem highly relates to the question whether 
$p\co \Acy_{2g} \to \Nil_{2g}$ is injective or not.
\end{remark}

\section{Refinements of the 
Johnson homomorphisms}
\label{ch:remark}
Finally, we summarize invariants of homology cylinders which extend 
the Johnson homomorphisms. For each $k \ge 2$, The $(k-1)^{\rm st}$ 
Johnson homomorphism for 
homology cylinders is nothing other than 
the composite $J_k \circ \sigma^{\mathrm{acy}} |_{\Hg[k]}$, 
where $\Hg[k]:= \Ker \left( \sigma_k\co  
\Hg \to \Aut_0 N_k \right)$, and we now denote it by 
$J_k$, for short. We can determine the image of $J_k$ by using 
\fullref{thm:gl} due to 
Garoufalidis--Levine and Habegger. See \cite{gl} for details. 

As seen in \fullref{ch:properties}, we have a refinement 
\[\widetilde{J}_k\co \Aut \Acy_{2g}[k] 
\to \Hom (F_{2g}, (\Gamma^k F_{2g})/(\Gamma^{2k-1} F_{2g}) )\] 
of the Johnson homomorphism $J_k$, so that 
we obtain an exact sequence
\[1 \longrightarrow \Hg[2k-1] \longrightarrow \Hg[k] 
\xrightarrow{\widetilde{J}_k \circ \sigma^{\mathrm{acy}}} 
\Hom (F_{2g}, (\Gamma^k F_{2g})/(\Gamma^{2k-1} F_{2g}) ).\]
On the other hand, we can construct maps which are essentially 
the same as $\widetilde{J}_k \circ \sigma^{\mathrm{acy}}$ 
by the following geometric method or group-homological one. 
This argument is based on Heap's idea in \cite{he} for 
the case of the mapping class group. 

The first one is obtained by applying 
the construction of the invariant $\theta$ in 
the previous section to the case of $\Hg[k]$. 
The same argument gives a homomorphism $\theta_k\co \Hg[k] \to H_3 (N_k)$ 
defined by
\[\theta_k ([M,i_+,i_-]) :=\widehat{f} ([C_{M}]) 
\quad \in H_3 (N_k),\]
and we can show that it is onto. 
Note that such a construction already appeared in \cite{gl}. 
\begin{thm}\label{bordism}
The kernel of $\theta_k$ is $\Hg[2k-1]$, 
namely we have an exact sequence
\[\begin{CD}
1 @>>> \Hg[2k-1] @>>> \Hg[k] @>\theta_k>> 
H_3 (N_k) @>>> 1.
\end{CD}\]
\end{thm}
\begin{proof}
In \cite{io}, 
Igusa--Orr showed that 
the homomorphism $H_3 (N_{2k-1}) \to H_3 (N_k)$ 
induced by the natural projection $N_{2k-1} \to N_k$ 
is trivial. From this, we see that 
$\Hg[2k-1] \subset \Ker \theta_k$. On the other hand, 
the induced homomorphism 
$\overline{\theta_k} \co \Hg[k]/\Hg[2k-1] \to H_3 (N_k)$ 
turns out to be a surjective one between 
free abelian groups of the same rank (see \cite{gl} and \cite{io}), 
which shows that 
$\overline{\theta_k}$ is an isomorphism. In particular, 
$\Hg[2k-1] = \Ker \theta_k$ follows. 
\end{proof}

The second one is obtained by generalizing Morita's 
refinement of the Johnson homomorphism for the mapping class group 
in \cite{mo} under the following modification. 
For $M \in \Hg[k]$, 
we write $\varphi_M := \sigma^{\mathrm{acy}}(M) \in 
\Aut_0 \Acy_{2g}$. 
Let $C_{\ast} (\Acy_{2g})$ be the chain complex 
of $\Acy_{2g}$ called the 
{\it normalized standard resolution} in \cite{br}. 
We take a 2--chain $c_0 \in C_2 (\Acy_{2g})$ satisfying 
$\partial c_0=-(\zeta)$ and whose image in $C_2 (\pi_1 
\Sigma_g)^{\mathrm{acy}}$ by the map from $\Sg$ to the closed surface 
$\Sigma_g$ 
coincides with the natural image of the fundamental cycle 
of $C_2 (\pi_1 \Sigma_g)$. For each $\varphi_M$, 
$c_0 - \varphi_M(c_0)$ gives a 2--cycle of $\Acy_{2g}$. 
Since $H_2 (\Acy_{2g})=0$, there exists a 3--chain 
$c_M \in C_3 (\Acy_{2g})$ satisfying 
$\partial c_M =c_0 - \varphi_M(c_0)$. 
Now $M \in \Hg[k]$, so that $\varphi_M$ acts on $N_k$ trivially. 
Hence if we write $\overline{c}_M$ for the image of $c_M$ 
in $C_3 (N_k)$, we obtain a 3--cycle $\overline{c}_M$ in $C_3 (N_k)$. 
We define $\widetilde{\tau}_k (M):=[\overline{c}_M] \in H_3 (N_k)$. 
A similar argument to one in \cite{mo} shows that 
$\widetilde{\tau}_k$ is well-defined and gives a refinement of $J_k$, 
where 
we need to use Igusa--Orr's result mentioned above 
instead of $H_3 (F_{2g})=0$. 

The same statement as \fullref{bordism} holds for 
$\widetilde{\tau}_k$. To show that $\Ker \widetilde{\tau}_k =\Hg[2k-1]$, 
we can apply Yokomizo's argument \cite{yo}, where the exact sequence 
\begin{align*}
H_3 (N_{k+1}) \to H_3 (N_k) \to H_1 (F_{2g}) \otimes 
(\Gamma^k F_{2g})/&(\Gamma^{k+1} F_{2g}) \to\\ 
&(\Gamma^{k+1} F_{2g})/(\Gamma^{k+2} F_{2g})\end{align*}
and Igusa--Orr's result are effectively used. The surjectivity of 
$\widetilde{\tau}_k$ follows from \fullref{thm:gl} and 
a direct sum decomposition of $H_3 (N_k)$ 
due to Igusa--Orr.

\bibliographystyle{gtart}
\bibliography{link}

\end{document}